\documentclass[reqno]{amsart}
\usepackage{amssymb}
\usepackage{eucal}
\usepackage{cite, url}
 \makeatletter

 \def\activeat#1{\csname @#1\endcsname}
 \def\def@#1{\expandafter\def\csname @#1\endcsname}
 {\catcode`\@=\active \gdef@{\activeat}}

\let\ssize\scriptstyle
\newdimen\ex@	\ex@.2326ex

 \def\requalfill{\cleaders\hbox{$\mkern-2mu\mathord=\mkern-2mu$}\hfill
  \mkern-6mu\mathord=$}
 \def\eqfill{$\m@th\mathord=\mkern-6mu\requalfill}
 \def\deffill{\hbox{$:=$}$\m@th\mkern-6mu\requalfill}
 \def\fiberbox{\hbox{$\vcenter{\hrule\hbox{\vrule\kern1ex
     \vbox{\kern1.2ex}\vrule}\hrule}$}}

 \newdimen\arrwd 
  \newdimen\minCDarrwd \minCDarrwd=2.5pc
   \setbox\z@\hbox{$\rightarrow\,$} \minCDarrwd=\wd\z@
 	
 \def\findarrwd#1#2#3{\arrwd=#3%
  \setbox\z@\hbox{$\ssize\;{#1}\;\;$}%
  \setbox\@ne\hbox{$\ssize\;{#2}\;\;$}%
  \ifdim\wd\z@>\arrwd \arrwd=\wd\z@\fi
  \ifdim\wd\@ne>\arrwd \arrwd=\wd\@ne\fi}
 \newdimen\arrowsp\arrowsp=0.375em  	
 \def\findCDarrwd#1#2{\findarrwd{#1}{#2}{\minCDarrwd}
    \advance\arrwd by 2\arrowsp}
 \newdimen\minarrwd 
 \setbox\z@\hbox{$\longrightarrow$} \minarrwd=\wd\z@

 \def\harrow#1#2#3#4{{\minarrwd=#1\minarrwd
   \findarrwd{#2}{#3}{\minarrwd}\kern\arrowsp
    \mathrel{\mathop{\hbox to\arrwd{#4}}\limits^{#2}_{#3}}\kern\arrowsp}}

 \def@]#1>#2>#3>{\harrow{#1}{#2}{#3}\rightarrowfill}
 \def@>#1>#2>{\harrow1{#1}{#2}\rightarrowfill}
 \def@<#1<#2<{\harrow1{#1}{#2}\leftarrowfill}
 \def@={\harrow1{}{}\eqfill}
 \def@:#1={\harrow1{}{}\deffill}
 \def@ N#1N#2N{\vCDarrow{#1}{#2}\UpDownarrow}
 \def\UpDownarrow{\uparrow\,\Big\downarrow}

\def\hookrightarrowfill{\hbox{$\lhook\joinrel$}\rightarrowfill}
\def@{hk>}#1>#2>{\harrow1{#1}{#2}\hookrightarrowfill}
\def@ c>#1>#2>{\harrow1{#1}{#2}\hookrightarrowfill}
\def\hookleftarrowfill{\leftarrowfill\hbox{$\joinrel\rhook$}}
\def@{hk<}#1<#2<{\harrow1{#1}{#2}\hookleftarrowfill}

 \def@.{\ifodd\row\relax\harrow1{}{}\hfill
   \else\vCDarrow{}{}.\fi}
 \def@|{\vCDarrow{}{}\Vert}
 \def@ V#1V#2V{\vCDarrow{#1}{#2}\downarrow}
 \def@ A#1A#2A{\vCDarrow{#1}{#2}\uparrow}
 \def@(#1){\arrwd=\csname col\the\col\endcsname\relax
   \hbox to 0pt{\hbox to \arrwd{\hss$\vcenter{\hbox{$#1$}}$\hss}\hss}}

 \def\squash#1{\setbox\z@=\hbox{$#1$}\finsm@@sh}
\def\finsm@@sh{\ifnum\row>1\ht\z@\z@\fi \dp\z@\z@ \box\z@}

 \newcount\row \newcount\col \newcount\numcol \newcount\arrspan
 \newdimen\vrtxhalfwd  \newbox\tempbox

 \def\innernewdimen{\alloc@1\dimen\dimendef\insc@unt}
 \def\measureinit{\col=1\vrtxhalfwd=0pt\arrspan=1\arrwd=0pt 
   \setbox\tempbox=\hbox\bgroup$}
 \def\setinit{\col=1\hbox\bgroup$\ifodd\row
   \kern\csname col1\endcsname
   \kern-\csname row\the\row col1\endcsname\fi}
 \def\findvrtxhalfsum{$\egroup
  \expandafter\innernewdimen\csname row\the\row col\the\col\endcsname
  \global\csname row\the\row col\the\col\endcsname=\vrtxhalfwd
  \vrtxhalfwd=0.5\wd\tempbox
  \global\advance\csname row\the\row col\the\col\endcsname by \vrtxhalfwd 
  \advance\arrwd by \csname row\the\row col\the\col\endcsname
  \divide\arrwd by \arrspan
  \loop\ifnum\col>\numcol \numcol=\col%
     \expandafter\innernewdimen \csname col\the\col\endcsname
     \global\csname col\the\col\endcsname=\arrwd
   \else \ifdim\arrwd >\csname col\the\col\endcsname
      \global\csname col\the\col\endcsname=\arrwd\fi\fi
   \advance\arrspan by -1 %
   \ifnum\arrspan>0 \repeat}
 \def\setCDarrow#1#2#3#4{\advance\col by 1 \arrspan=#1 
    \arrwd= -\csname row\the\row col\the\col\endcsname\relax
    \loop\advance\arrwd by \csname col\the\col\endcsname
     \ifnum\arrspan>1 \advance\col by 1 \advance\arrspan by -1%
     \repeat
    \squash{\mathop{
     \hbox to\arrwd{\kern\arrowsp#4\kern\arrowsp}}\limits^{#2}_{#3}}}
 \def\measureCDarrow#1#2#3#4{\findvrtxhalfsum\advance\col by 1%
   \arrspan=#1\findCDarrwd{#2}{#3}%
    \setbox\tempbox=\hbox\bgroup$}
 \def\vCDarrow#1#2#3{\kern\csname col\the\col\endcsname
    \hbox to 0pt{\hss$\vcenter{\llap{$\ssize#1$}}%
     \Big#3\vcenter{\rlap{$\ssize#2$}}$\hss}\advance\col by 1}

 \def\setCD{\def\harrow{\setCDarrow}%
  \def\\{$\egroup\advance\row by 1\setinit}
  \m@th\lineskip3\ex@\lineskiplimit3\ex@ \row=1\setinit}
 \def\endsetCD{$\egroup}
 \def\dr@p#1\\{\findvrtxhalfsum\advance\row by 2 \measureinit}
 \def\measure{\bgroup
  \def\harrow{\measureCDarrow}%
  \def\\##1{\ifx##1\endmeasure\endmeasure\else\expandafter\dr@p\fi}%
  \row=1\numcol=0\measureinit}
 \def\endmeasure{\findvrtxhalfsum\egroup}

 \newcount\savedcount	
 \def\LCD#1\end{\savedcount=\count11
   \measure#1\endmeasure
   \vcenter{\setCD#1\endsetCD\kern\medskipamount}%
   \global\count11=\savedcount\end}

 \newenvironment{CD}{\let\at=@\catcode`\@=\active\LCD}{\catcode`\@=12\relax}

 \makeatother
 \usepackage{fancyhdr} 

\def\thetitle{The Canonical Model of a Singular Curve}
\newcommand{\emdash}{\unskip\penalty10000\thinspace
        ---\penalty-500\thinspace\ignorespaces}

\def\UThin{\penalty\@M \thinspace\ignorespaces}
\def\(#1){{\let~=\UThin\rm(#1)}}
\def\tsum{\textstyle\sum}

\DeclareMathOperator{\Proj}{Proj}
\DeclareMathOperator{\Spec}{Spec}
\DeclareMathOperator{\Supp}{Supp}
\DeclareMathOperator{\Torsion}{{\it Torsion}}
\DeclareMathOperator{\cok}{{\it Cok}}

\DeclareMathOperator{\rIm}{{\rm Im}}
 \DeclareMathOperator{\sHom}{{\it Hom}}
 \DeclareMathOperator{\Hom}{Hom}
\DeclareMathOperator{\Ext}{Ext}
\DeclareMathOperator{\Sym}{Sym}
\DeclareMathOperator{\sSym}{{\it Sym}}

\let\?=\overline
\let\:=\colon
\let\dg=\dagger
\let\bb=\mathbb

\let\into=\hookrightarrow
\let\mc=\mathcal
\let\mf=\mathfrak
\let\onto=\twoheadrightarrow
\let\ox=\otimes

\let\ve=\varepsilon

\let\x=\times
\let\xto=\xrightarrow
\let\wh=\widehat

\def\risom{\buildrel\sim\over{\smashedlongrightarrow}}
 \def\smashedlongrightarrow{\setbox0=\hbox{$\longrightarrow$}\ht0=1.25pt\box0}
\newcommand{\fm}{\mf m}   
\newcommand{\sM}{\mc M}   
\newcommand{\sC}{\mc C}   

\theoremstyle{plain}
 \newtheorem{thm}{Theorem}[section]
 \newtheorem{cor}[thm]{Corollary}
 \newtheorem{lem}[thm]{Lemma}
 \newtheorem{prp}[thm]{Proposition}

\theoremstyle{definition}
 \newtheorem{dfn}[thm]{Definition}
 \newtheorem{eg}[thm]{Example}

  \newtheorem{rmk}[thm]{Remark}

\makeatletter
 \@addtoreset{equation}{thm}
\makeatother

\def\mylistparam
 {\setbox0=\hbox{\rm(a)\kern0.5em}
  \setlength{\labelwidth}{\wd0}
  \setlength{\topsep}{\smallskipamount}%
  \setlength{\itemsep}{0pt}\setlength{\parsep}{0pt}%
  \setlength{\itemindent}{0pt}%
  \setlength{\leftmargin}{\parindent}%
  \addtolength{\leftmargin}{\wd0}%
 }%
\renewenvironment{enumerate}%
  {\begin{list}{\rm(\alph{enumi})}%
  {\usecounter{enumi}\mylistparam}%
  }%
 {\end{list}}

\begin{document}

\title\thetitle

\author[S. L. Kleiman]{Steven Lawrence Kleiman}
 \address
 {Department of Mathematics, MIT \\
 77 Mass. Ave.\\
 Cambridge, MA 02139, USA}
 \email{kleiman@math.mit.edu}

\author[R. V. Martins]{Renato Vidal Martins}
\address{Departamento de Matem\'atica, ICEx, UFMG
Av. Ant\^onio Carlos 6627,
30123-970 Belo Horizonte MG, Brazil}
\email{renato@mat.ufmg.br}
\thanks{The second author was supported in part by CNPq
 grant number PDE 200999/2005-2.}
\date{\today}



\begin{abstract}
We give refined statements and modern proofs of Rosenlicht's results
about the canonical model $C'$ of an arbitrary complete integral curve
$C$.  Notably, we prove that $C$ and $C'$ are birationally equivalent if
and only if $C$ is nonhyperelliptic, and that, if $C$ is
nonhyperelliptic, then $C'$ is equal to the blowup of $C$ with respect
to the canonical sheaf $\omega$.  We also prove some new results: we
determine just when $C'$ is rational normal, arithmetically normal,
projectively normal, and linearly normal.
\end{abstract}

\maketitle

\vspace{-\smallskipamount}  

\section{Introduction}\label{intro}
Let $C$ be a complete integral curve of arithmetic genus $g\ge2$ defined
over an algebraically closed field of arbitrary characteristic.  Its
canonical model $C'$ was introduced by Rosenlicht at the end of his
paper \cite{R} (based on his 1950 Harvard thesis under Zariski) where he
introduced the dualizing sheaf $\omega$.  Here, we give modern proofs of
Rosenlicht's results about $C'$; also, we determine just when $C'$ is
rational normal, arithmetically normal, projectively normal, and
linearly normal.

Rosenlicht constructed the {\it canonical model\/} $C'$ as follows.  He
\cite[p.\,188 top]{R} proved that $H^0(\omega)$ defines a
base-point-free linear series on the normalization $\?C$ of $C$.  He
formed the corresponding map $\?\kappa\: \?C\to \bb P^{g-1}$, and took
its image to be $C'$.

Rosenlicht \cite[p.\,188 top]{R} called $C$ ``quasihyperelliptic'' if
$\?\kappa$ is not birational onto $C'$.  He \cite[Thm.\,15, p.\,188]{R}
proved that $C$ is quasihyperelliptic iff there is some map
$\lambda\:C\to \bb P^1$ of degree 2.  Nowadays, it is more
common to call $C$ {\it
  hyperelliptic\/} if such a map $\lambda$ exists; so Rosenlicht's
result is just our Proposition~\ref{prBir}.  Furthermore, it is implicit
in Rosenlicht's work, and it is easy to prove, see St\"ohr's discussion
\cite[p.\,96, top]{St99} or our Proposition~\ref{prHyp2}, that if
$\lambda$ exists, then it is unique.  In fact, then $\lambda$ is induced
by the canonical map $\?\kappa\: \?C\to \bb P^{g-1}$, and its image $C'$
is equal to the rational normal curve $N_{g-1}$ of degree $g-1$; furthermore,
$\omega\simeq\lambda^*\mc O_{\bb P^1}(g-1)$.

Suppose $C$ is hyperelliptic.  Then $\omega\simeq\lambda^*\mc O_{\bb
  P^1}(g-1)$; whence, $\omega$ is invertible, so $C$ is Gorenstein.
Rosenlicht \cite[p.\,188 top]{R} reasoned in essentially this way, as
did St\"ohr \cite[p.\,96, top]{St99} and as do we in proving
Proposition~\ref{prHyp2}.  On the other hand, Homma \cite[Cor.\ 3.3,
  p.\,31]{Ho99} reproved that $C$ is Gorenstein, but he proceeded
differently; he obtained and used an explicit equation for a plane model
of $C$.

As to a nonhyperelliptic $C$, first Rosenlicht \cite[Cor.\ and Thm.\,17,
p.\,189]{R} proved these three statements: (1) if $C$ is
nonhyperelliptic and Gorenstein, then $\?\kappa$ induces an isomorphism
$\kappa\: C\risom C'$; (2) furthermore, then $C'$ is {\it extremal\/};
that is, its genus is maximal for its degree, which is $2g-2$; and (3)
conversely, every extremal curve of degree $2g-2$ in $\bb P^{g-1}$ is
nonhyperelliptic and Gorenstein of genus $g$, and is its own canonical
model.  Rosenlicht's proofs involve relating global invariants.  We
give similar proofs of these statements in our Theorem~\ref{thNonh}.

Unaware of Rosenlicht's work, several authors have reproved various form
of (1).  The first proofs were given by Deligne and Mumford in 1969, by
Sakai in 1977, and by Catanese in 1982, according to
Catanese~\cite[p.\,51]{Ca82}.  Their work was motivated by the study of
families of curves, and they allowed $C$ to be reducible, but required
it to be connected in a strong sense, and to have only mild Gorenstein
singularities.

In 1973, Mumford and Saint-Donat proved (1) for a smooth $C$, using the
Jacobian of $C$.  In 1983, Fujita \cite[p.\,39]{Fu83} asserted their
proof works virtually without change for any Gorenstein $C$.  Then
Fujita \cite[Thm.\,(A1), p.\,39]{Fu83} gave another proof for any
Gorenstein $C$, involving ideas from Mumford's version of Castelnuovo
Theory.

In Remark~\ref{rm3rdpf}, we explain a variant of the latter proof,
involving ideas from the version of Castelnuovo Theory developed by
Arbarello et al.\ \cite[pp.\,114--117]{ACGH}.  In this way, both Fujita,
in his Theorem (A1), and we, in our Proposition~\ref{prGpn}, obtain
more: namely, $C'$ is {\it projectively normal\/}; that is, for every
$n\ge1$, the hypersurfaces of degree $n$ cut out a complete linear
series.  In particular, $C'$ is {\it linearly normal\/}; that is, the
hyperplanes cut out a complete series.  In fact, the converse holds in
general; indeed, using Castelnuovo Theory, we prove Lemma~\ref{leLnPn},
which asserts that, whether $C$ is Gorenstein or not, $C'$ is linearly
normal iff it is projectively normal.

In 1986, Hartshorne \cite[Thm.\,1.6, p.\,379]{Ha86} gave yet another
proof of (1); he showed that the complete linear series of canonical
divisors of $C$ ``separates points and tangent vectors.''  We give a
somewhat similar proof in Section~\ref{sc_iso}, and obtain a stronger
statement, Theorem~\ref{prGiso}, which  is virtually
in Rosenlicht's paper; namely,  $\?\kappa\: \?C\to
\bb P^{g-1}$ induces an open embedding $\kappa\:G\into C'$ where $G$ is
the Gorenstein locus of $C$, the largest open set on which $\omega$ is
invertible.

 Rosenlicht's last result \cite[Thm.\,17, p.\,189]{R} is his main
theorem about $C'$; it  asserts that, if
$C$ is nonhyperelliptic, then the birational map between $C$ and $C'$ is
regular on $C'$.  In fact, Rosenlicht's proof nearly yields a more
refined result, which is our Theorem~\ref{thRMT}.  It asserts that
$\?\kappa$ induces an isomorphism $\wh\kappa\:\wh C\risom C'$ where $\wh
C$ is the blowup of $C$ with respect to $\omega$ in the sense of
Definition~\ref{dfBlwp}.  Instead, Rosenlicht \cite[p.\,191, last
  line]{R} worked with the subsheaf $\omega'\subset \omega$ generated by
$H^0(\omega)$, but apparently, he was unaware that $\omega'=\omega$, a
fundamental discovery made by Eisenbud, Harris, Koh, and Stillman
\cite[p.\,536, mid]{EKS}.

Rosenlicht's proof of his main theorem involves some hard local algebra,
which reduces the general statement to (1) above.  We give a similar
proof.  We can have no purely local proof until we find a purely local
condition equivalent to the nonhyperellipticity of $C$.  A sufficient
condition is that $C$ have a point of multiplicity at least 3 by our
Proposition~\ref{prHyp2}(1), but this multiplicity condition is hardly
necessary.

Rosenlicht's proof of (1)--(3) involves his version of Clifford's
Theorem \cite[Thm.\,16, p.\,188]{R}, the first version for a singular
curve.  It concerns an invertible sheaf $\mc F$ on $C$ such that
$h^0(\mc F)\ge1$ and $h^1(\mc F)\ge1$; it asserts {\it Clifford's
Inequality}
\begin{equation*}\label{eq}
h^0(\mc F) + h^1(\mc F) \le g+1,
\end{equation*}
and it describes when equality holds.

Correspondingly, we prove a version of Clifford's Theorem,
Theorem~~\ref{leClif}.  It is more general, as we prove the above
inequality for any torsion-free sheaf $\mc F$ of rank~1 on $C$.  The
added generality is due to Kempf~\cite[pp.\,25, 32]{Ke71}, and we
present his proof, which is short and has not fully appeared in print
before.

If equality holds, then $H^0(\mc F)$ generates $\mc F$.  This result is
due to Eisenbud et al.\ \cite[p.\,536, mid]{EKS}.  They derived it in
a few lines from the above inequality for the subsheaf $\mc F'\subset
\mc F$ generated by $H^0(\mc F)$.  We reproduce their proof, as the
result is fundamental.  St\"ohr \cite[Thm.\,3.2, p.\,123]{S}
rediscovered the case $\mc F = \omega$; his proof is  different,
and we discuss it in Remark~\ref{rmSt}.

If equality holds and if $\mc F$ is invertible, but not isomorphic to
either $\mc O_C$ or $\omega$, then $C$ is hyperelliptic.  This statement
was proved by Rosenlicht, and we reprove it by modifying the standard
proof in the case where $C$ is smooth \cite[pp.\,344--345]{Ha86}.  In
this case, $\mc F$ is isomorphic to the pullback of $\mc O(n)$ under the
map $\lambda\:C\to \bb P^1$ with $n:=h^0(\mc F)-1$.  Conversely, if $C$
is hyperelliptic and if $0\le n\le g-1$, then for $\mc F:=\lambda^*\mc
O(n)$, equality holds in Clifford's Inequality.

Surprisingly, equality can hold in Clifford's Inequality, yet $\mc F$ is
neither invertible nor isomorphic to $\omega$.  Cases were discovered
and classified by Eisenbud et al.\ \cite[Thm.\,A(c), p.\,533]{EKS}.
Namely, $C$ is rational, and $\mc F$ is isomorphic to the sheaf
generated over $\mc O_C$ by $H^0(\bb P^1, \mc O(n))$ inside the pushout
of $\mc O(n)$ under the normalization map $\nu\:\bb P^1\to C$ with
$n:=h^0(\mc F)-1$.  Moreover, $C$ is not hyperelliptic; in fact, the
canonical map $\?\kappa\:\bb P^1\to \bb P^{g-1}$ is the Veronese
embedding, so that $C$ and $C'$ are birational, and $C'$ is the rational
normal curve $N_{g-1}$ of degree $g-1$.  Furthermore, $C$ is {\it nearly
  normal;} that is, $C$ has a unique multiple point $P$, and its maximal
ideal sheaf $\sM_{\{P\}}$ is equal to the conductor $\sC$.

In addition, Eisenbud et al.\ \cite[Rmk.\ p.\,533]{EKS} observed that
this $C$ is isomorphic to a curve of degree $2g+1$ in $\bb P^{g+1}$ that
lies on the cone $S$ over the rational normal curve $N_g$ of degree $g$
in $\bb P^g$ and that has a unique multiple point at the vertex;
moreover, the canonical map corresponds to the projection from a ruling.
Conversely, if $C$ is isomorphic to a curve of degree $2g+1$ on $S$ with
a unique multiple point at the vertex, then $C$ is as described in the
preceding paragraph; this converse was discovered by the second
author~\cite[Thm.\,2.1, p.\,461]{M}, and we reprove it differently and
in a stronger form as part of Theorem~\ref{coRnc}.

Thus we obtain three characterizations of a $C$ whose canonical model
$C'$ is the rational normal curve $N_{g-1}$: (1) a $C$ with an $\mc F$,
other than $\mc O_C$ or $\omega$, 
for which equality holds in Clifford's Inequality; (2) a $C$ isomorphic
to a curve of degree $2n+1$ in $\bb P^{n+1}$ lying on the cone $S$ over
the rational normal curve $N_n$ of degree $n$ in $\bb P^n$ for some
$n\ge2$; a posteriori, $n=g$; and (3) a $C$ that is either hyperelliptic
or else rational and nearly normal.

If $C$ is nearly normal with unique multiple point $P$, then the local
ring of $C$ at $P$ is of an interesting sort, which was introduced and
studied by Barucci and Fr\"oberg.  Namely, they
\cite[p.\,418]{BF} termed a 1-dimensional local Cohen--Macaulay ring
with finite integral closure {\it almost Gorenstein} if its
Cohen--Macaulay type satisfies a certain relation, recalled below in
Definition~\ref{dfnlyG}.  So, in Definition~\ref{dfnlyG}, we term $C$
{\it nearly Gorenstein} if the non-Gorenstein locus $C-G$ consists of a
single point whose local ring is almost Gorenstein.  Theorem~\ref{thAn}
asserts, in particular, that if $C$ is nearly normal, but
non-Gorenstein, then it is nearly Gorenstein.

More generally, Theorem~\ref{thAn} characterizes a non-Gorenstein $C$,
rational or not, whose canonical model $C'$ is {\it arithmetically
normal\/}; that is, its homogeneous coordinate ring is normal.  Namely,
if $C$ is non-Gorenstein, then these seven conditions are equivalent:
(a) $C'$ is arithmetically normal; (b) $C'$ is smooth and projectively
normal; (c) $C'$ is smooth and linearly normal; (d) $C'$ is smooth and
extremal; (e) $C'$ is of degree $g+\?g-1$ where $\?g$ is the genus of
$\?C$, the normalization; (f) $C$ is nearly normal; and (g) $C$ is
nearly Gorenstein, and $\smash{\wh C}$ is smooth, where $\smash{\wh C}$
is the blowup with respect to $\omega$.  Furthermore, if these
conditions hold, then, at its unique multiple point, $C$ is of
multiplicity $g-\?g+1$ and of embedding dimension $g-\?g+1$; thus $C$ is
of {\it maximal embedding dimension,} as the embedding is always bounded
by the multiplicity according to Lipman's Corollary
1.10~\cite[p.\,657]{Li71}.  In order to prove our theorem, we use
Castelnuovo Theory and some propositions due to Barucci and Fr\"oberg
\cite{BF}.

Our final result, Theorem~\ref{thNmlty} applies to even more $C$.
Namely, it asserts that, if $C$ is non-Gorenstein, then these six
conditions are equivalent: (a) $C'$ is projectively normal; (b) $C'$ is
linearly normal; (c) $C'$ is extremal; (d) $C'$ is of degree $g+g'-1$
where $g'$ is the genus of $C'$; (e) $C$ is nearly Gorenstein; and (f)
$C'=\Spec\bigl(\sHom(\sM_{\{P\}},\,\sM_{\{P\}})\bigr)$ where $\sM_{\{P\}}$ is the
maximal ideal sheaf of some point $P$ off the Gorenstein locus.
Furthermore, if (f) holds, then $C$ is of maximal embedding dimension at
$P$ iff $C'$ is Gorenstein.  In order to prove this theorem, we use
Rosenlicht's Main Theorem, that if $C$ is nonhyperelliptic, then there is a
canonical isomorphism $\wh\kappa\:\wh C\risom C'$.

In short, Section 2 develops the preliminary theory of the canonical
model $C'$, including the basic theory of hyperellipticity and some
results about the degree $d'$ of $C'$.  Section~3 proves Clifford's
theorem, and applies it to characterize the case where $C'$ is the
rational normal curve $N_{g-1}$.  Section~4 characterizes the
nonhyperelliptic and Gorenstein $C$, and proves that, for an arbitrary
nonhyperelliptic $C$, the canonical map induces an open embedding
$\kappa\:G\into C'$, where $G$ is the Gorenstein locus of $C$.
Section~5 develops Castelnuovo Theory, and applies it to characterize
the case where $C'$ is arithmetically normal.  Finally, Section~6 proves
Rosenlicht's Main Theorem, and applies it to characterize the case where
$C'$ is projectively normal, or equivalently, linearly normal.

\section{The canonical model}\label{sc_cm}
\label{sec:cm}

Let $C$ be an arbitrary complete integral curve over an algebraically
closed base field $k$ of arbitrary characteristic.  Let $g$ denote its
arithmetic genus, and assume $g\ge2$.  Let $\omega_C$, or simply
$\omega$, denote the canonical sheaf (dualizing sheaf).

Recall that $H^0(\omega)$ generates $\omega$.  In full generality, this
fundamental result was discovered and proved by Eisenbud, Harris, Koh,
and Stillman \cite[p.\,536 mid]{EKS} in 1988.  Their argument is
recalled in the proof of Lemma~\ref{leClif}; in fact, with $\mc
F:=\omega$, the lemma yields the result.  The result was rediscovered,
by St\"ohr \cite[Thm.\,3.2, p.\,123]{S} in 1993, and proved in a
different way, which is described in Remark~\ref{rmSt}.

In the special case that $C$ is Gorenstein (that is, $\omega$ is
invertible), this result was obtained by a number of authors.  The first
was Rosenlicht \cite[p.\,187 bot]{R}; in fact, he proved only that
$H^0(\omega)$ generates the pullback of $\omega$ to the normalization of
$C$; however, it follows immediately, via Nakayama's lemma, that
$H^0(\omega)$ generates $\omega$, because $\omega$ is invertible.
Catanese \cite[Thm.\,D, p.\,75]{Ca82} rediscovered the result in 1982; in
fact, he worked with a reducible $C$, and found conditions guaranteeing
that $H^0(\omega)$ generates $\omega$. 

Fujita \cite[Thm.\,(A1), p.\,39]{Fu83} rediscovered the Gorenstein case
in 1983.  In fact, he claimed that Mumford and Saint-Donat's 1973 proof
\cite[Prp.\,(1.5), p.\,160]{Sa73} works virtually without change,
although they assumed $C$ to be smooth.  Then Fujita gave his own proof.
Finally, Hartshorne \cite[Thm.\,1.6, p.\,379]{Ha86} gave yet another
proof in 1986; furthermore, he \cite[Rmk.\,1.6.2, p.\,380]{Ha86} cited
Fujita's work and Catanese's work.

\begin{dfn}\label{dftor}
As a matter of notation, given any integral scheme $A$ and any coherent
sheaf $\mc F$ on $A$, let $\Torsion(\mc F)$ denote its torsion subsheaf.
And given any map $\alpha\:A\to C$ and any sheaf $\mc G$ on $C$, set
	$$\mc O_A\mc G := \alpha^*\mc G\big/\!\Torsion(\alpha^*\mc G).$$
\end{dfn}

\begin{dfn}\label{dfCM}
Let $\nu\: \?C \to C$ denote the normalization map.  Then $\mc
O_{\?C}\omega$ is invertible, and is generated by $H^0(\omega)$; hence,
there is a natural nondegenerate map
 $$\?\kappa\: \?C\to \bb P^{g-1}.$$

 If $C$ is Gorenstein, then $\omega$ is invertible and
 generated by $H^0(\omega)$.  So then $\?\kappa$ factors uniquely:
 $$\?\kappa = \kappa\nu \text{ where }\kappa\: C\to \bb P^{g-1}.$$
 
  Call the above maps $\?\kappa$ and $\kappa$ the {\it canonical maps\/}
  of $C$.  Call their common image the {\it canonical model} of $C$, and
  denote it by $C'$.  Set $d':=\deg C'$. 

When appropriate, let $\?\kappa$ and $\kappa$ also denote the induced
maps:
 $$\?\kappa\: \?C\to C'\text{ and }\kappa\:C\to C'$$
 Furthermore, denote the arithmetic genus of $C'$ by $g'$, and that of
 $\?C$ by $\?g$.
\end{dfn}

\begin{rmk}\label{rmCM}
Under the conditions of Definition~\ref{dfCM}, if conversely $\?\kappa =
\kappa\nu$, then $\omega$ and $\kappa^*\mc
 O_{C'}(1)$ are equal, because both are equal to the subsheaf of
 $\nu_*\mc O_{\?C}\omega$ generated by $H^0(\omega)$; whence, then $C$
 is Gorenstein.
\end{rmk}

\begin{dfn}\label{dfHyp} As usual, call $C$ {\it hyperelliptic\/} if
there is some map of degree 2 $$\lambda\:C\to \bb P^1.$$ Otherwise, call
$C$ {\it nonhyperelliptic}.  \end{dfn}

\begin{eg}\label{egPQ}
  Suppose $C$ is a plane quartic.  Then $\omega =\mc O_C(1)$.  Hence
  $C'=C$ and $\kappa=1_C$.  Furthermore, $C$ is nonhyperelliptic by
  Proposition~\ref{prHyp2}\,(1) below.

  Suppose also that $C$ is 3-nodal.  Then $\?C=\bb P^1$ and $\mc
  O_{\?C}\omega=O_{\?C}(4)$.  Yet $C$ depends on three moduli.  Thus the
  position of $H^0(\omega)$ in $H^0(O_{\?C}\omega)$ is crucial for
  $\?\kappa$ and $C'$.
\end{eg}

\begin{prp}\label{prHyp2} Assume $C$ is hyperelliptic.

\(1) Then there is an isomorphism $\omega\simeq\lambda^*\mc O_{\bb
  P^1}(g-1)$.

 \(2) Then $C$ is Gorenstein with double points at worst, and
  $\deg\kappa=\deg\?\kappa=2$.

 \(3) Then $\kappa=\ve\lambda$ where $\ve\:\bb
   P^1\to\bb P^{g-1}$ is isomorphic to the Veronese embedding.

   \(4) Then $\lambda$ is uniquely determined, up to an automorphism of
   $\bb P^1$.
\end{prp}
\begin{proof}
  Given a (closed) point $P\in C$, let $u$ be a uniforming parameter at
  $\lambda(P)$ on $\bb P^1$.  Then $\dim(\mc O_P/\langle u\rangle)\le2$
  since $\deg\lambda=2$.  Hence $P$ is of multiplicity at most 2.
   
  Set $\mc L:=\lambda^*\mc O_{\bb P^1}(g-1)$.  Then (a) $\deg\mc L=2g-2$
  because $\deg\lambda=2$, and (b) $ h^0(\mc L)\ge g$ since the
  natural map $\mc O_{\bb P^1}(g-1)\to\lambda_*\mc L$ is plainly
  injective.  Now, (a) implies $\chi(\mc
  L)=g-1$.  So (b) implies $h^1(\mc L)\ge1$; whence, by duality, there
  is a nonzero map $w\:\mc L\to \omega$.  Since $C$ is integral, $w$ is
  injective; whence, $\cok(w)=0$ since $\chi(\mc L) =\chi(\omega)$.
  Thus $w$ is bijective; so (1) holds.

  The remaining assertions hold
  just because $\omega\simeq\lambda^*\mc O_{\bb P^1}(g-1)$ and
  $\deg\lambda=2$.
 \end{proof}

\begin{dfn}\label{dfCond}
 As a further matter of notation, set
 $$\mc O:=\mc O_C\text{ and }  \?{\mc O}:= \nu_*\mc O_{\?C}.$$
 Let $\sC$ denote the conductor of  $\?{\mc O}$ into $\mc O$.
Given a point $P\in C$, set
 $$\delta_P := \dim (\?{\mc O}_P/\mc O_P)
 \text{ and } \eta_p := \delta_P - \dim (\mc O_P/\sC_P).$$
 Furthermore, set
 $$\delta := \tsum_{P\in C}\delta_p = h^0(\?{\mc O}/\mc O) \text{ and }
 \eta := \tsum_{P\in C}\eta_P = h^0(\mc O/\sC).$$
 Finally, set
 $$\?\omega:=\nu_*\omega_{\?C}\text{ and }
 \?{\mc O}\omega := \nu_*(\mc O_{\?C}\omega).$$
\end{dfn}

 In the next lemma, the main assertion is the equation $\sC\omega =
 \?\omega$.  It was proved implicitly by Rosenlicht
 \cite[pp.\,177--180 bot]{R} and Serre \cite[\S\,11, p.\,80]{Sr}, and
was proved explicitly by St\"ohr \cite[Prp.\,2.2, p.\,113]{S}.  It is
 proved here a bit differently.

\begin{lem}\label{leCshf} We have $\sC\omega=\sC\?{\mc O}\omega =
  \?\omega\subset \omega\subset\?{\mc O}\omega$.  Furthermore, given
  $P\in C$, there is an $x\in\omega_P$ such that $\?{\mc O}_Px = (\?{\mc
    O}\omega)_P$; any such $x$ satisfies $\sC_Px=\?\omega_P$.
 \end{lem}
 \begin{proof}
   Plainly $\omega\subset\?{\mc O}\omega$. Now, $\?\omega=\sHom(\?{\mc
     O},\,\omega)$ by general principles \cite[Ex.\,7.2(a), p.\,249]{H};
   so  $\sC\omega\subset \?\omega \subset \omega$.  Given
   $P\in C$ and $y\in\?\omega_P$, we have to prove $y\in \sC_P\omega_P$.

   Since $\?{\mc O}_P$ is a semilocal Dedekind domain, it's a UFD; so
there is an $x\in\omega_P$ such that $\?{\mc O}_Px = (\?{\mc O}\omega)_P$.  Fix
such an $x$.  Then $y=ax$ for some $a\in\?{\mc O}_P$.  We have to prove
$a\in\sC_P$, for then $y\in\sC_Px\subset\sC_P\omega_P$.  So, given
$b\in\?{\mc O}_P$, we have to prove $ab\in\mc O_P$.

By general principles, $\mc O\risom\sHom(\omega,\,\omega)$; indeed, the
natural map is injective, whence bijective since source and target have
the same Euler characteristic by duality.  So, given $z\in \omega_P$, we
have $abz\in(\?{\mc O}\omega)_P$ since $(\?{\mc O}\omega)_P$ is an
$\?{\mc O}_P$-module, and we have to prove $abz\in\omega_P$.

Say $z=cx$ where $c\in\?{\mc O}_P$.  Now,  $y\in\?\omega_P$, and
$\?\omega_P$ is an $\?{\mc O}_P$-module.  So $bcy\in\?\omega_P$.  But
$bcy=abcx=abz$.  Thus $abz\in\?\omega_P\subset\omega_P$, as desired.    
\end{proof}

In the next lemma, the first assertion is well known.  It was proved by
Rosenlicht \cite[Thm.\,10, p.\,179]{R} first, and his proof was repeated
by Serre \cite[\S\,11, p.\,80]{Sr}.  The proof is repeated here,
because, with one additional line, it yields Formula~(\ref{eqEta2}).
Alternatively, this formula holds because, as observed by Eisenbud et
al.\ \cite[p.\,535 mid]{EKS}, the residue map induces a perfect pairing
on $\mc O_P/\sC_P\x (\?{\mc O}\omega)_P/\omega_P$.

\begin{lem}\label{leEta}
  Fix $P\in C$.  Then $\eta_P\ge0$, with equality iff $\mc O_P$ is
Gorenstein.  Also,
\begin{equation}\label{eqEta2}
  \eta_p = \delta_P - \dim ((\?{\mc O}\omega)_P/\omega_P).
 \end{equation}
 \end{lem}
 \begin{proof}
   By Lemma~\ref{leCshf}, there is an $x\in\omega_P$ such that $\sC_Px
   =\?\omega_P$.  The latter equation is plainly equivalent to the
   injectivity of following map:
  $$\mc O_P/\sC_P\to\omega_P/\?\omega_P
  \text{ defined by } f\mapsto fx.$$ 
  The image is $\mc O_Px/\?\omega_P$.  By duality,
  $\dim(\omega_P/\?\omega_P) = \delta$.  Hence,
\begin{equation}\label{eqEta1}
  \eta_P = \dim ( \omega_P/\mc O_Px).
\end{equation}
 
Therefore, $\eta_P\ge0$, and if equality holds, then $\mc O_P$ is
Gorenstein.  Conversely, if $\mc O_P$ is Gorenstein, then there is a
$y\in\omega_P$ such that ${\mc O}_Py= \omega_P$, and so $\?{\mc O}_Py =
(\?{\mc O}\omega)_P$; whence, by Lemma~\ref{leCshf}, we may take $x:=y$,
and so $\eta=0$.

Finally, $\dim((\?{\mc O}\omega)_P/\mc O_Px) =\delta_P$ as $(\?{\mc
  O}\omega)_P = \?{\mc O}_Px$.  So (\ref{eqEta1}) yields
(\ref{eqEta2}).
  \end{proof}

  The following lemma is essentially Eisenbud et al.'s \cite[Lem.\,2,
  p.\,534]{EKS}, and their proof is essentially the alternative proof
  here.

\begin{lem}\label{prDeg}
 We have $\deg\mc O_{\?C}\omega = 2g-2-\eta$. 
\end{lem}
\begin{proof}
  The Riemann--Roch Theorem and the birational invariance of
  $\chi(\bullet)$ yield
$$ \deg{\mc O}_{\?C}\omega = \chi({\mc O}_{\?C}\omega) -\chi(\mc O_{\?C})
    =\chi(\?{\mc O}\omega)-\chi(\?{\mc O}).$$
But Formula~(\ref{eqEta2}) and Definition~\ref{dfCond} yield
$\chi(\?{\mc O}\omega/\omega)=\delta-\eta$ and $\chi(\?{\mc O}/\mc
O)=\delta$.  So the additivity of $\chi(\bullet)$
  and the duality equation $\chi(\omega)= -\chi(\mc O)$ yield the assertion:
$$\deg{\mc O}_{\?C}\omega= (\chi(\omega)+\delta-\eta) - (\chi(\mc O)
  +\delta ) =  2g-2-\eta.$$

 Alternatively, for each $P\in C$, there is an $x_P\in\omega_P$ such that
 $\?{\mc O}_Px_P = (\?{\mc O}\omega)_P$ by Lemma~\ref{leCshf}.
 Plainly, the various $\mc O_Px_P$ are the stalks of an invertible sheaf
 $\mc L\subset \omega$ such that $\nu^*\mc L\risom \mc O_{\?C}\omega$.
 Then $\deg\mc L = \deg\omega -\sum_{P\in C}\dim(\omega_P/\mc O_Px_P)$,
 and also, $ \deg\mc L =  \deg\nu^*\mc L= \deg\mc O_{\?C}\omega$;
whence, (\ref{eqEta1}) yields the asserted formula.
 \end{proof}

\begin{lem}\label{leDeg} We have $\deg\mc O_{\?C}\omega \le 2g-2$, with
equality iff $C$ is Gorenstein. 
 \end{lem}
\begin{proof}
 The assertion follows immediately from Lemma~\ref{prDeg} and
Lemma~\ref{leEta}.
  \end{proof}

\begin{lem}\label{leDeg2}
  We have $d' \le (2g-2)/(\deg\?\kappa)$, with equality iff $C$ is
  Gorenstein.
\end{lem}
\begin{proof}
  As $d' = (\deg\mc O_{\?C}\omega) / (\deg\?\kappa)$,
Lemma~\ref{leDeg} yields the assertion.
\end{proof}

\begin{prp}\label{prBir}
  If $\deg\?\kappa=1$, then $C$ is nonhyperelliptic, and conversely.
\end{prp}
\begin{proof}
  If $\deg\?\kappa=1$, then $C$ is nonhyperelliptic by
  Proposition~\ref{prHyp2}(1). Conversely, suppose $\deg\?\kappa\ge2$.
  Then $d'\le g-1$ by Lemma~\ref{leDeg2}.  But $C'$ spans $\bb
  P^{g-1}$.  Hence $C'$ is the rational normal curve of degree $g-1$ by
  a well-known old theorem \cite[p.\,18]{B}.  So
  Lemma~\ref{leDeg2} implies $\deg\?\kappa=2$ and $C$ is Gorenstein.
  Hence $\kappa\:C\to C'$ exists and is of degree 2.  So $C$ is
  hyperelliptic with $\lambda:=\kappa$ since $C'\simeq\bb P^1$.
\end{proof}

\begin{prp}\label{prdegC'}
If $C$ is nonhyperelliptic, then $d' = 2g-2-\eta$. 
\end{prp}
\begin{proof}
  As $d' = (\deg\mc O_{\?C}\omega) / (\deg\?\kappa)$,
Proposition~\ref{prBir} and Lemma~\ref{prDeg} yield the assertion.
\end{proof}

\begin{dfn}\label{dfnnl} Call $C$ {\it nearly normal\/} if $h^0(\mc
O/\sC)=1$, that is, if $C$ is has a unique multiple point
$P$ and its maximal ideal sheaf $\sM_{\{P\}}$ is equal to the conductor
$\sC$.  \end{dfn}

\begin{prp}\label{prLbd}
  Suppose that $C$ is nonhyperelliptic.  If $C$ is smooth, then $\?g=g$
  and $d'= 2g-2$.  If $C$ is singular, then $d'\ge g+\?g-1$,
  with equality iff $C$ is nearly normal.
 \end{prp}
\begin{proof}
 Proposition \ref{prBir} yields $\deg\?\kappa=1$.  Hence $d'= \deg
 \mc O_{\?C}\omega$.  So Lemma \ref{prDeg} yields $d'= 2g-2-\eta$.
 But clearly, $\eta = \delta -h^0(\mc O/\sC)$ and
 $\delta =
 g-\?g$.  Hence
\begin{equation*}\label{eqDC'}
	d' = g+\?g-2+h^0(\mc O/\sC).
\end{equation*}
If $C$ is smooth, then $\?C=C$; whence, $\?g=g$ and  $h^0(\mc
O/\sC)=0$. If $C$ is singular, then  $h^0(\mc
O/\sC)\ge1$, with equality iff  $C$ is nearly normal.
 \end{proof}

\section{Rational normal models}\label{sc_rat}
Preserve the setup introduced at the beginning of Section~\ref{sc_cm}
and after Proposition~\ref{prHyp2}.  In this section, the main result is
Theorem~\ref{coRnc}, which characterizes the case in which the canonical
model $C'$ is equal to the rational normal curve.

Rosenlicht \cite[Thm.\,16, p.\,188]{R} was the first to prove the next
lemma, Clifford's theorem, for any $C$, but for an invertible $\mc F$.
For an arbitrary $\mc F$, Kempf \cite[p.\,32]{Ke71} was the first to
prove the bound in (1).  Eisenbud, Harris, Koh, and Stillman
\cite[Thm.\,A, p.\,532]{EKS} used Kempf's argument (credited via
\cite[p.\,544]{E89}) and Eisenbud's study of determinantal varieties
\cite[Thm.\,B, p.\,537]{EKS} to prove (1)--(5) and to generalize (5) to
an arbitrary $\mc F$, thereby discovering a surprising new case; more
about this case is said in Remark~\ref{rm2}.  Here, (5) is proved for
any $C$ by adapting the proof in \cite[pp.\,344--345]{Ha86}, or what is
virtually the same, that in \cite[\S1]{Sa73}.

\begin{lem}[Clifford's Theorem]\label{leClif}
  Let $\mc F$ be torsion-free sheaf of rank $1$ on $C$ such that
  $h^0(\mc F)\ge1$ and $h^1(\mc F)\ge1$.

\(1) Then $h^0(\mc F) + h^1(\mc F) \le g+1$.
 If equality holds, then $H^0(\mc F)$ generates $\mc  F$.

\(2)  Equality holds in {\rm(1)} and $h^0(\mc F)=1$ iff $\mc F\simeq
  \mc O_C$.

\(3)  Equality holds in {\rm(1)} and $h^1(\mc F)=1$ iff  $\mc F\simeq
  \mc \omega$.

\(4) Assume $C$ is hyperelliptic.  Then equality holds in {\rm(1)} iff
  $\mc F\simeq\lambda^*\mc O_{\bb P^1}(n)$ with $0\le n\le g-1$; if so, then
  $h^0(\mc F)= n+1$.

\(5) Assume equality holds in {\rm(1)} and $\mc F$ is invertible.
  Assume  either $h^0(\mc F)=2$ or else $h^0(\mc F)\ge 3$
  and $h^1(\mc F)\ge 2$.  Then $C$ is hyperelliptic.
\end{lem}
\begin{proof}
  In (1), let us prove the bound following Kempf~\cite[pp.\,25, 32]{Ke71}.  Observe
  that the pairing
\begin{equation*}\label{eqCl1}
H^0(\mc F) \x \Hom(\mc F,\,\omega) \to H^0(\omega)
\end{equation*}
is nondegenerate; that is, if $(f,u)\mapsto 0$, then $f=0$ or $u=0$.
But, given any three $k$-vector spaces $A$, $B$, and $C$ of dimensions
$a$, $b$, and $c$ and given any nondegenerate pairing $A\x B\to C$, then
$a+b\le c+1$; indeed, $A\x B$ may be viewed canonically as the set of
$k$-points on a cone in the affine space whose $k$-points are $A\ox B$,
and this cone meets, only at the origin, the affine space whose
$k$-points form the kernel of the induced map $A\ox B\to C$.  (This part
of Kempf's proof has not appeared in print before.  However, according
to Arbarello et al.\ \cite[p.\,135]{ACGH}, the bound $a+b\le c+1$ itself
was ``used'' by H. Hopf in 1940/41.)  Thus  the bound  holds.

  In (1), assume equality holds.  Following Eisenbud et
  al. \cite[p.\,536, mid]{EKS}, form the subsheaf $\mc G\subset\mc F$
  generated by $H^0(\mc F)$.  Consider the induced sequence
\begin{equation*}\label{eqCl2}
 H^0(\mc G) \xto{u} H^0(\mc F)\to H^0(\mc F/\mc G)\to H^1(\mc G)
 \xto{v} H^1(\mc F)\to H^1(\mc F/\mc G).
\end{equation*}
By construction, $u$ is an isomorphism; so $h^0(\mc G)=h^0(\mc F)$.
Now, $h^0(\mc F)\ge1$ and $\mc F$ is torsion-free sheaf of rank $1$;
hence, $\mc F/\mc G$ has finite support.  So $H^1(\mc F/\mc G)=0$.
Hence $v$ is surjective, and so $h^1(\mc G)\ge h^1(\mc F)$.  Hence
\begin{equation*}\label{eqCL3}
h^0(\mc G)+h^1(\mc G)\ge h^0(\mc F)+h^1(\mc F).
\end{equation*}
The left side is at most $g+1$ by the bound with $\mc G$
for $\mc F$.  The right side is equal to $g+1$ by assumption.  Hence
$h^1(\mc G)= h^1(\mc F)$.  Hence $v$ is an isomorphism.  Therefore,
$H^0(\mc F/\mc G)=0$.  Hence $\mc G=\mc F$.  Thus $H^0(\mc F)$ generates
$\mc F$.

To prove (2) and (3), assume equality holds in (1).  If $h^0(\mc F)=1$,
then there is a nonzero map $\mc O_C\to \mc F$; it is plainly injective,
so bijective as $\chi(\mc O_C)=\chi(\mc F)$.  If $h^1(\mc F)=1$, then
similarly, there is a bijection $\mc F\risom\omega$. 
The converses plainly hold.

To prove (4), assume $C$ is hyperelliptic.  First, also assume equality
holds in (1).  Set $n:=h^0(\mc F)-1$.  Then $h^1(\mc F)=g-n$.  So, by
hypothesis, $g-n\ge1$.  So $\mc O_{\bb P^1}(g-n-1)$ has a section that
vanishes nowhere on the image under $\lambda$ of the singular locus
$\Sigma$ of $C$.  Let $D$ be the section's divisor of zeros, and form
the corresponding sequence
\begin{equation}\label{eqCL6}
0\to\mc O_{\bb P^1}\to\mc O_{\bb P^1}(g-n-1)\to\mc O_D\to0.
\end{equation}
Set  $\mc G:=\mc F\ox \lambda^*\mc O_{\bb P^1}(g-n-1)$.
Pull (\ref{eqCL6}) back to $C$, and then tensor with  $\mc F$.  The
result is an exact sequence
\begin{equation*}\label{eqCL7}
0\to\mc F\to\mc G\to\mc O_{\lambda^{-1}D}\to0,
\end{equation*}
because $\mc F$ is invertible along $\lambda^{-1}D$ since
$\lambda(\Sigma)\cap D=\emptyset$.  Hence
\begin{equation}\label{eqCL8}
\chi(\mc G) = \chi(\mc F) + \deg(\lambda)\cdot\deg(D)
 = n+1-(g-n)+2(g-n-1) = g-1.
\end{equation}

Consider the natural pairing
\begin{equation*}\label{eqCL4}
H^0(\mc F)\x H^0(\mc O_{\bb P^1}(g-n) \to H^0(\mc G).
\end{equation*}
Plainly, it is nondegenerate.  So  $h^0(\mc G)\ge (n+1)+(g-n)-1=g$.  Hence
$h^1(\mc G)\ge 1$ owing to (\ref{eqCL8}).  So
\begin{equation*}\label{eqCL5}
h^0(\mc G)+h^1(\mc G) \ge g+1.
\end{equation*}
 But the opposite inequality holds by (1) for $\mc G$.  Hence equality
holds, and $h^1(\mc G)=1$.  Therefore, $\mc G\simeq \omega$ by (3).  But
$\omega=\lambda^*\mc O_{\bb P^1}(g-1)$ by Proposition~\ref{prHyp2}(2).
Thus $\mc F\simeq\lambda^*\mc O_{\bb P^1}(n)$, and plainly $0\le n\le
g-1$.

Conversely, also assume $\mc F\simeq\lambda^*\mc O_{\bb P^1}(n)$ with $0\le
n\le g-1$.  Then
\begin{equation*}\label{eqCL51}
h^0(\mc F)\ge h^0(\mc O_{\bb P^1}(n))=n+1.
\end{equation*}
Now, $\omega\simeq\lambda^*\mc O_{\bb P^1}(g-1)$ by
Proposition~\ref{prHyp2}(2).  So $\sHom(\mc
F,\,\omega)\simeq\lambda^*\mc O_{\bb P^1}(g-1-n)$.  Therefore, $\Hom(\mc
F,\,\omega)\supset H^0(\mc O_{\bb P^1}(g-1-n))$; whence, $h^1(\mc F)\ge
g-n$.  So
\begin{equation*}\label{eqCL61}
h^0(\mc F) + h^1(\mc F) \ge g+1,
\end{equation*}
But the opposite inequality holds by (1).  So equality holds in (1).
Thus (4) holds.

To prove (5), assume equality holds in (1) and $\mc F$ is invertible.
Then $H^0(\mc F)$ generates $\mc F$ by (1).  Set $n:=h^0(\mc F) -1$.
Then there is a map $\rho\:C \to \bb P^n$ such that $\mc F=\rho^*\mc
O_{\bb P^n}(1)$ and $\rho C$ spans $\bb P^n$.  Moreover, the equality in
(1) yields $h^1(\mc F)= g-n$.

First, also assume $h^0(\mc F)= 2$.  Then $n=1$ and $h^1(\mc F)= g-1$.
So the Riemann--Roch theorem yields $\deg\mc F =2$.  So $\rho\:C \to \bb
P^1$ has degree 2.  Thus $C$ is hyperelliptic.

Finally, also assume $h^0(\mc F)\ge 3$ and $h^1(\mc F)\ge 2$.  Let
$\Sigma$ be the singular locus of $C$.  Fix a point $P$ of $C$ off
$\rho^{-1}\rho\Sigma$.  Since $h^1(\mc F)\ge 2$, there is a nonzero map
$u\:\mc F\to \omega$ that vanishes at $P$.  Now, $n\ge2$ and $\rho C$
spans $\bb P^n$; so $\rho C$ is not a line.  So there is a point $Q$ of
$C$ such that (1) $\rho Q$ does not lie on the line through $\rho P$ and
$\rho S$ for any $S\in \Sigma$ and such that (2) $u$ is not bijective at
$Q$.  Then the line through $\rho P$ and $\rho Q$ does not contain $\rho
S$ for any $S\in \Sigma$.  So there is a hyperplane $H$ containing $\rho
P$ and $\rho Q$, but not $\rho S$ for any $S\in \Sigma$.

Set $D:=\rho^{-1}H$.  Then $D$ is an effective divisor, which contains
$P$ and $Q$, but no point of $\Sigma$; furthermore, $\mc O_C(D)=\mc F$.
Fix a map $f\:\mc O_C\to\mc F$ whose locus of zeros is $D$.  Then $f$
induces a map $h\:\sHom(\mc F,\,\omega)\to \omega$.  Form the map
\begin{equation*}\label{eqCL9}\textstyle
(u,-h)\:\mc F \bigoplus \sHom(\mc F,\,\omega)\to \omega.
\end{equation*}
Denote its kernel and image by $\mc F'$ and $\mc G$, and form the short
exact sequence
\begin{equation}\label{eqCL10}\textstyle
0\to\mc F'\to\mc F \bigoplus \sHom(\mc F,\,\omega)\to \mc G\to 0.
\end{equation}
Plainly, $\mc F'$ and $\mc G$ are  torsion-free sheaf of rank $1$.

Let $R\in \Sigma$ be arbitrary.  Then $R\notin D$.  So $f\:\mc O_C\to\mc
F$ is bijective at $R$.  Hence $h\:\!\sHom(\mc F,\,\omega)\to \omega$ is
bijective at $R$.  Therefore, $\mc F'\to\mc F$ is bijective at $R$.  It
follows that $\mc F'$ is invertible and that $\mc F'\to\mc F$ is
injective.  Furthermore, $\mc F'\to\mc F$ is not bijective at $Q$
because $u$ is not bijective at $Q$.  Therefore, $h^1(\mc F')\ge h^1(\mc
F)\ge2$ and $\deg\mc F'<\deg\mc F$.

The map $u\:\mc F\to \omega$ induces a map $v\:\mc O_C\to\sHom(\mc
F,\,\omega)$.  With it, form the map
\begin{equation*}\label{eqCL11}\textstyle
  \binom fv\:\mc O_C\to\mc F \bigoplus \sHom(\mc F,\,\omega).
\end{equation*}
It is injective as $f$ is (or as $v$ is).  Plainly, the composition
$(u,-h)\binom fv$ vanishes.  Hence $\binom fv$ factors through $\mc F'$,
so induces a section of $\mc F'$.  This section is nonzero, and it
vanishes at $P$ since both $f$ and $v$ do.  Hence $\mc F'$ is
nontrivial and $h^0(\mc F')\ge1$.

On the category of torsion-free $\mc O_C$-modules $\mc M$, the functor
$\sHom(\mc M,\,\omega)$ is dualizing; see, for example, from \cite[Thm.
21.21]{E94}.  Now, by definition, $h = \sHom(f,\,\omega)$ and
$v=\sHom(u,\,\omega)$.  Hence there are canonical isomorphisms
$f=\sHom(h,\,\omega)$ and $u=\sHom(v,\,\omega)$.
Therefore, the functor turns the map $\binom fv$ into the map
\begin{equation*}\label{eqCL12}\textstyle
(u,h)\:\mc F \bigoplus \sHom(\mc F,\,\omega)\to \omega.
\end{equation*}
Plainly, the image of $(u,h)$ is equal to that of
$(u,-h)$, namely, $\mc G$. 

Since $\binom fv$ factors through $\mc F'$, dually $(u,h)$ factors
through $\sHom(\mc F',\,\omega)$.  The map into $\sHom(\mc F',\,\omega)$
is surjective, because the dual of Sequence~(\ref{eqCL10}) is exact, as
$\mc G$ is torsion free.  The map out of $\sHom(\mc F',\,\omega)$ is
injective, as this sheaf is torsion free.  Hence $\sHom(\mc
F',\,\omega)$ is equal to the image of $(u,h)$, so to $\mc G$.
Therefore, $h^1(\mc F')=h^0(\mc G)$.

Since the sequence~(\ref{eqCL10}) is exact, we therefore have
\begin{equation*}\label{eqeqCL13}
h^0(\mc F) + h^1(\mc F) \le h^0(\mc F') + h^1(\mc F').
\end{equation*}
By hypothesis, the left-hand side is equal to $g+1$.  By (1), the
right-hand side is at most $g+1$.  Hence the right-hand side is equal to
$g+1$.  But $\mc F'$ is nontrivial and $h^0(\mc F')\ge1$.  Hence (2)
implies $h^0(\mc F')\ge2$. 
 But $h^1(\mc F')\ge2$ and $\mc F'$ is invertible.
 Thus the hypotheses of (5) are fulfilled
with $\mc F'$ in place of $\mc F$.  But $\deg\mc F'<\deg\mc F$.
Therefore, by induction on $\deg\mc F$, we may conclude that $C$ is
hyperelliptic.  Thus (5) holds.
\end{proof}

\begin{prp}\label{prg=2} Assume $g=2$.  Then $C$ is hyperelliptic iff
$C$ is Gorenstein.
  \end{prp}
 \begin{proof}
   If $C$ is Gorenstein, then Clifford's Theorem \ref{leClif}(5) with
   $\mc F:=\omega$ implies that $C$ is hyperelliptic.  The converse
   holds by Proposition~\ref{prHyp2}(1).
 \end{proof}

\begin{lem}\label{leCone}
  Let $S\subset \bb P^{n+1}$ be the cone, say with vertex $P$, over the
  rational normal curve $N_n$ of degree $n$ in $\bb P^n$ with $n\ge2$.
  Let $D$ be a
  curve on
  $S$ of degree $2n+1$.  Then $D$ contains $P$, and has arithmetic
  genus $n$; furthermore, the canonical map of $D$ is isomorphic to the
  projection from a ruling of $S$.  If $P$ is simple on $D$, then $D$ is
  hyperelliptic.  If $P$ is multiple, then its multiplicity is $n+1$, it
  is resolved by blowing up, and it is the only multiple point of $D$;
  moreover, then the maximal ideal sheaf of $P$  is equal to the
  conductor of $D$, and $D$ is rational, but not Gorenstein.
\end{lem}
\begin{proof}
  Let $\wh S$ be the blowup of $S$ at $P$, and $E$ the exceptional
  divisor.  Then $\wh S$ is a rational ruled surface, and $E$ is a
  section with $E^2=-n$, according to \cite[Ch.\,V, Sec.\,2]{H}.  Let
  $f$ be a ruling.  Let $H$ be the pullback of a hyperplane section.
  Then $H\equiv E+nf$; indeed, $H\cdot f=1$ as $f$ maps isomorphically
onto a line,
  and $H\cdot E=0$ as $E$ contracts to $P$.  Let $\wh D$ be the strict
  transform of $D$.  Say $\wh D\equiv aE+bf$.  Then $b=2n+1$ since $\wh
  D\cdot H=\deg D$.  But $b\ge an$ by \cite[Cor.\,2.18(b)]{H}; hence,
  $a=1,2$.

  Suppose $a=2$.  Then $\wh D\cdot E=1$; hence, the map $\wh D\to D$ is
  an isomorphism, and $P\in D$ is simple.  And $\wh D\cdot f=2$; hence,
  the projection of $\wh S$ onto the base $\bb P^1$ induces a map $\wh
  D\to\bb P^1$ of degree 2.  So $D$ is hyperelliptic, with the
  projection from $P$ as the map $\lambda\:D\to \bb P^1$ of degree 2.
  But the projection of $D$ from a ruling is equal to $\lambda$ followed
  by the projection of the rational normal curve $N_n$ 
  from a point; hence, the  projection from a ruling  is isomorphic to the
  canonical map by Proposition~\ref{prHyp2}(2).  Furthermore, the
  arithmetic genus of $D$ is equal to $1+(D+K)\cdot D/2$ where $K$ is a
  canonical divisor on $\wh S$.  But $K\equiv -2E-(2+n)f$ by
  \cite[Cor.\,{\bf V},\,2.11]{H}.  Therefore, $D$ is of genus $n$.

Finally, suppose $a=1$.  Then $\wh D\cdot E=n+1$; hence, the
multiplicity of $P$ on $D$ is $n+1$.  And $\wh D\cdot f=1$; hence, the
projection of $\wh S$ onto the base $\bb P^1$ induces an isomorphism
$\wh D\risom \bb P^1$, and so $P$ is resolved by blowing up, $P$ is the
only multiple point of $D$, and $D$ is rational.

Denote by $A$ the local ring of $D$ at $P$, by $M$ its maximal ideal,
and by $\?A$ its normalization.  Then $\dim(\?A/M\?A)$ is equal to the
multiplicity of $A$, so to $n+1$.  Denote by $n^\dg$ the arithmetic
genus of $D$.  Then $\dim(\?A/M)=n^\dg+1$ since $D$ is rational and
$P$ is its only multiple point.  Hence $n^\dg\ge n$, with equality iff
$M=M\?A$.

Set $h_i:=h^i(\mc O_{D}(1))$.  Then $h_0\ge n+2$ because $D$ lies in no
hyperplane $H$ of $\bb P^{n+1}$; else, $D$ lies in $H\cap S$, so is a
union of lines.  Hence the Riemann--Roch Theorem yields the inequality:
$n+2\le(2n+1)+1-n^\dg+h_1$, or $n^\dg\le n+h_1$.  If $h_1>0$, then
Clifford's Inequality, Lemma~\ref{leClif}(1), yields the opposite
inequality: $n+2+h_1\le n^\dg+1$, or $n + h_1 <n^\dg$.  Thus $h_1=0$ and
so $n^\dg\le n$.

Combined, the above two paragraphs yield the equations: $n^\dg=n$ and
$M=M\?A$.  Therefore, $M$ is an ideal in $\?A$, so is the conductor of
$\?A$ into $A$.  As $P$ is the only multiple point, the maximal ideal
sheaf of $P$ is equal to the conductor of $D$.

The definition of $\eta$ now yields $\eta=n-1$.  But $n\ge2$.  So
$\eta\ge1$.  Hence $D$ is not Gorenstein by Lemma~\ref{leEta}.
Therefore, $D$ is nonhyperelliptic by Proposition~\ref{prHyp2}(1).
Denote by $D'$ the canonical model of $D$.  Then $\deg D'= n-1$ by
Proposition~\ref{prLbd} since $D$ is rational and nonhyperelliptic.  But
$D'$ spans $\bb P^{n-1}$.  Hence $D'$ is the rational normal curve
$N_{n-1}$ by a well-known old theorem \cite[p.\,18]{B}.  Since the
projection of $D$ from a ruling of $S$ is also birational onto the the
rational normal curve $N_{n-1}$, the projection is isomorphic to the
canonical map.  \end{proof}

\begin{thm}\label{coRnc}
  The following three conditions are equivalent:
\begin{enumerate}
 \item either $C$ is hyperelliptic, or $C$ is rational and nearly normal;

 \item $C$ is isomorphic to a curve of degree $2n+1$ in $\bb P^{n+1}$
that lies on the cone $S$ over the rational normal curve $N_n$ of degree
$n$ in $\bb P^n$ for some $n$ at least $2$;

 \item $C'$ is equal to the rational normal curve $N_{g-1}$ of degree $g-1$ in
   $\bb P^{g-1}$.
 \end{enumerate} If \(b) holds, then  $n=g$, and the
canonical map corresponds to the projection from a ruling.  \end{thm}
\begin{proof}

Condition (b) implies (a) and (c) and the final assertion by
Lemma~\ref{leCone}.

Assume that (c) holds and that $C$ is nonhyperelliptic.  Then
$\deg\?\kappa=1$ by Lemma~\ref{prBir}.  Hence $\?\kappa\:\?C\to C'$ is
an isomorphism since $C'$ is smooth.  So $\?g=0$.  But $d' = g-1$.
So (a) holds by Proposition~\ref{prLbd}.  Thus (c) implies (a).

Finally, assume (a), and let's prove (b).  First, suppose $C$ is
hyperelliptic, and let $\lambda\:C\to \bb P^1$ be the map of degree 2.
Pick a simple point $P\in C$, and set $\mc L:= \lambda^*\mc O_{ \bb
P^1}(g)\ox \mc O_C(P)$.  Then $\deg(\mc L)=2g+1$.  So $\mc L$ is very
ample and $h^1(\mc L)=0$ by Lemma~\ref{leAmp}.  The Riemann--Roch
Theorem yields $h^0(\mc L)=g+2$.  So $\mc L$ provides an embedding
$\gamma\:C\into \bb P^{g+1}$ with $\deg(\gamma(C))= 2g+1$.  Since $\mc
L:= \lambda^*\mc O_{ \bb P^1}(g)\ox \mc O_C(P)$, projection from
$\gamma(P)$ is the map provided by $\lambda^*\mc O_{ \bb P^1}(g)$.  So
its image is the rational normal curve $N_g$.  Hence $\gamma(C)$ lies on
the cone over $N_g$ with vertex at $\gamma(P)$.

Second, suppose $C$ is rational, $C$ has a unique multiple point $P$,
and $\sM_{\{P\}}=\sC$.  Then $\?C\simeq \bb P^1$.  Take a coordinate
function $x$ on $\?C$ that is finite on $\nu^{-1}P$.  Set $\?A:=k[x]$.
Let $A$ be the ring of $\nu(\Spec(\?A))\subset C$.  Form the conductor
of $\?A$ into $A$, say it's the principal ideal $f\?A$.  Then $f\?A$ is
equal to the maximal ideal of $P$ in $A$ since $\sM_{\{P\}}=\sC$.  Hence
$A=k+f\?A$.  Further, $\deg(f)=g+1$ since $\dim(\?A/A)=g$ as $C$ is
rational of genus $g$.  Say $f=x^{g+1}+a_1x^g+\dotsb$.  Then
$x^{g+1}=f-a_1x^g-\dotsb$.  It follows that $A$ is generated as a
$k$-algebra by $f,\,xf,\dotsc,\,x^{g}f$.

Consider the map $\?\gamma\:\?C\to\bb P^{g+1}$ given by
$$\?\gamma(x):=(f,\,xf,\dotsc,\,x^{g}f,\,1).$$ Then $\?\gamma$ factors
through a map $\gamma\:C\to\bb P^{g+1}$ because $x^if\in A$.  Further,
$\gamma$ is an embedding on $\Spec(A)$ because $A$ is generated as a
$k$-algebra by $f,\,xf,\dotsc,\,x^{g}f$, and $\gamma$ is an embedding at
infinity because $f/xf=1/x$.  Clearly, $\gamma(P)=(0,\dotsc,0,1)$ and
projection from $\gamma(P)$ maps $\gamma(C)$ birationally onto $N_g$.
Hence $\gamma(C)$ lies on the cone over $N_g$ with vertex at
$\gamma(P)$.  Thus (b) holds with $n:=g$.
 \end{proof}

\begin{rmk}\label{rm2}
 Part of Theorem~\ref{coRnc} was known already.  Eisenbud et al.\
\cite[Rmk.\ p.\,533]{EKS} noted that, if $C$ is rational, if $C$ has a
unique multiple point $P$, and if $\sM_{\{P\}}=\sC$, then $C$ is
isomorphic to a curve of degree $2g+1$ in $\bb P^{g+1}$ that lies on the
cone $S$ over the rational normal curve $N_g$ of degree $g$ in $\bb P^g$
and $C$ is multiple at the vertex.  The converse was discovered by
the second author~\cite[Thm.\,2.1, p.\,461]{M}, who sketched an elementary
computational proof.

Interest stems from Clifford's Theorem.  Indeed, Eisenbud et
al.\  proved that these $C$ are just the curves that possess a
noninvertible torsion-free sheaf $\mc F$ of rank $1$ such that $h^0(\mc
F)\ge2$ and $h^1(\mc F)\ge2$ and $h^0(\mc F)+h^1(\mc F)=g+1$.  Further,
as in the hyperelliptic case, for each integer $n$ with $1\le n\le g-1$,
there is a unique such $\mc F$ with $h^0(\mc F)= n+1$; namely, in terms of the
normalization map $\nu\: \?C \to C$, the $\mc O_C$-module $\mc F$ is
isomorphic to the submodule of $\nu_*\mc O_{\?C}(n)$ generated by the
vector space $H^0(\mc O_{\?C}(n))$.  \end{rmk}

\section{The Gorenstein locus}\label{sc_iso}
 Preserve the general setup introduced at the beginning of
Section~\ref{sc_cm} and after Proposition~\ref{prHyp2}.  In this
section, the main result is Theorem~\ref{prGiso}, which asserts that, if
$C$ is nonhyperelliptic, then its canonical map $\?\kappa: \?C\to \bb
P^{g-1}$ induces an open embedding of its Gorenstein locus $G\subset C$
into its canonical model $C'$.  The proof involves the blowup $\wh C$ of
$C$ with respect to $\omega$, which is introduced in
Definition~\ref{dfBlwp}.

First, Theorem~\ref{thNonh} treats the special case in which $C$ is
nonhyperelliptic and Gorenstein; in this case,  $\?\kappa$ induces an
isomorphism, $\kappa\:C\risom C'$.  This result was proved by Rosenlicht
\cite[Thm.\,17, p.\,189]{R}, who used the Riemann--Roch Theorem and
Clifford's Theorem on $C'$ much as here.  In addition,
Theorem~\ref{thNonh} gives three necessary and sufficient numerical
conditions for this case to occur.

Furthermore, Theorem~\ref{thNonh} asserts that, in this case, $C'$ is
{\it extremal\/}; that is, its genus is maximal for its degree, which is
$2g-2$.  The theorem also asserts the converse: every
extremal curve of degree $2g-2$ in $\bb P^{g-1}$ is Gorenstein of genus
$g$, and is its own canonical model.  These statements too were
discovered by Rosenlicht \cite[Cor., p.\,189]{R}, and proved much as
here.

Neither Theorem~\ref{thNonh} nor Theorem~\ref{prGiso} depends logically
on the other, and their proofs are rather different in nature; the
first is global, the second local.

\begin{lem}\label{leextr}
Assume that there is a nondegenerate embedding of $C$ in $\bb P^r$ with
degree $d$.  Set $h_i:=h^i(\mc O_C(1))$.  If $d<2r$, then $g\le d-r$ and
$h_1=0$.  If $d=2r$, then either $h_1=0$ and $g\le r$, or else $h_1=1$
and $g=r+1=h_0$ and $\mc O_C(1)\simeq \omega$.
\end{lem}
\begin{proof}
Since $C\subset\bb P^r$ is nondegenerate, $r\le h_0-1$.  By the
Riemann--Roch Theorem, $h_0-1 = d-g+h_1$.  Hence $g\le d-r+h_1$.  And so $g<
r+h_1$ if $d<2r$.

If $h_1>0$, then Clifford's Theorem~\ref{leClif}(1) yields $h_0+h_1\le
g+1$.  But $r\le h_0-1$.  So if $h_1>0$, then $r+h_1\le g$.  Thus if
$d<2r$, then $h_1=0$ and so $g\le d-r$.
 
Assume therefore $d=2r$.  Then $r\le h_0-1 = 2r-g+h_1$.  So $g\le
r+h_1$.  Thus $g\le r$ if $h_1=0$.  Assume therefore $h_1\ge1$ also.
Then $r+h_1\le g$.  Hence $g=r+h_1$.  So $h_0-1=2r-(r+h_1)+h_1=r$ and
$h_0+h_1=g+1$.

Thus, if $h_1=1$, then $g=r+1=h_0$, and by Lemma \ref{leClif}(3), then $\mc
O_C(1))\simeq \omega$.

Finally, if $h_1\ge2$, then the proof of Lemma \ref{leClif}(5) implies
$H^0(\mc O_C(1))$ defines a map $C\to \bb P^r$ of degree 2, contrary to
hypothesis.  Thus this case does not occur.
\end{proof}

\begin{dfn}\label{dfextr}
A nondegenerate (reduced and irreducible) curve in $\bb P^r$ is said to
be {\it extremal\/} if its (arithmetic) genus is maximal among all curves
of its degree.
\end{dfn}

\begin{thm}\label{thNonh} The following six conditions are equivalent:
\begin{enumerate}\label{enNonh}
\item  $C$ is nonhyperelliptic and Gorenstein;
\item  $d' = 2g-2$;
\item  $g'=g$.
\item  $d' = g'+g-2$;
\item  $\kappa\:C\to C'$ exists, and is an isomorphism.
\item  $C$ is isomorphic to a curve $C^\dg\subset \bb P^r$ that is
nondegenerate, of degree $2r$, and extremal, for some $r\ge2$.
 \end{enumerate}
 If \(f) holds, then the isomorphism is equal to $\kappa$; in
particular,  $C^\dg=C'$ and $r=g-1$. 
\end{thm}
\begin{proof}
  Lemma~\ref{leDeg} implies $d' = 2g-2$ iff $\deg\?\kappa=1$ and
  $C$ is Gorenstein.  So Proposition~\ref{prBir} implies (a) and (b)
  are equivalent.

  Set $h_i:=h^i(\mc O_{C'}(1))$.  Then $g\le h_0$ because $C'$ is
  nondegenerate in $\bb P^{g-1}$ by construction.

Suppose (a) and (b) hold.  Apply the Riemann--Roch Theorem to $\mc
O_{C'}(1)$, getting
$$g \le (2g-2) + (1-g') + h_1, \text{ or } g'+1 \le g+h_1.$$
Now, (a) implies $\kappa\:C\to C'$ exists, and by Proposition
~\ref{prBir}, is birational; whence, $g\le g'$.  Hence $h_1\ge 1$.
So $h_0+h_1\le g'+1$ by Clifford's Theorem, \ref{leClif}(1).  Thus
$$ g'+1 \le g+h_1\le h_0+h_1\le g'+1;$$
whence, $g=h_0$ and $h_0+h_1= g'+1$.  But $g\ge2$, so $h_0\ge2$.  So, if
$h_1\ge2$, then $C'$ is hyperelliptic by Clifford's Theorem,
\ref{leClif}(5); hence, $C$ is hyperelliptic too, contrary to (a).  Hence
$h_1=1$.  Thus (c) holds.

Conversely, suppose (c) holds.  By hypothesis, $g\ge2$.  So $g'\neq0$.
Hence $C$ is nonhyperelliptic by Proposition~\ref{prHyp2}(3).  So
$\deg\?\kappa = 1$ by Proposition~\ref{prBir}.  So $d'\le 2g-2$,
with equality iff $C$ is Gorenstein by Lemma~\ref{leDeg}.  Apply the
Riemann--Roch Theorem to $\mc O_{C'}(1)$: since $g\le h_0$, and since
$g'=g$ by (c), we get
$$g \le (2g-2) + (1-g) + h_1, \text{ or } 1 \le h_1,$$
with equality only if $C$ is Gorenstein.  Hence, $g+h_1\le g'+1$ by
Clifford's Theorem, \ref{leClif}(1).  But $g'=g$.  Hence $h_1=1$, and so
$C$ is Gorenstein.  Thus (a) holds.

Plainly, (b) and (c) imply (d).  Conversely, suppose (d) holds.  Then
$C$ is nonhyperelliptic; otherwise, Proposition~\ref{prHyp2}(3) implies
$g'=0$ and $d' = g-1$, contradicting (d).  So $\deg\?\kappa = 1$ by
Proposition~\ref{prBir}.  So $d'\le 2g-2$, with equality iff $C$ is
Gorenstein by Lemma~\ref{leDeg}; whence, (d) implies $g'\le g$, with
equality iff $C$ is Gorenstein.  Apply the Riemann--Roch Theorem to $\mc
O_{C'}(1)$: since $g\le h_0$, we get
 $$g \le (g'+g-2) + (1-g') + h_1, \text{ or } 1 \le h_1.$$
 Hence $g+h_1\le g'+1$ by Clifford's Theorem \ref{leClif}(1), and so
 $g\le g'$.  Hence $g=g'$, and so $C$ is Gorenstein.  Thus (a) and (c)
 hold.

Suppose (a) and (c) hold.  Then, as noted above, $\kappa\:C\to C'$
exists and is birational.  But $g'=g$ by (c); hence, $\kappa$ is an
isomorphism.  Thus (e) holds.  Conversely, (e) trivially implies (c).
Thus  (a)--(e) are equivalent.

Suppose (a)--(e) hold, and set $r:=g-1$.  Then $C$ is isomorphic to
$C'\subset \bb P^r$ by (e).  Also, $C'$ is nondegenerate by
construction, and $C'$ is of degree $2r$ by (b).  If $r=1$, then $C'$ is
equal to $\bb P^1$, and so $C'$ is of degree 1; hence, $r\ge2$.
Furthermore, $C'$ is of genus $r+1$ by (c); hence, $C'$ is extremal,
because any nondegenerate curve of degree $2r$ in $\bb P^r$ is of
genus at most $r+1$ by Lemma~\ref{leextr}.  Thus (f) holds.

Finally, suppose (f) holds; that is, $C^\dg\subset \bb P^r$ is extremal
of degree $2r$.  But, we just proved that, given any Gorenstein curve of
genus $r+1$, its canonical model is an extremal curve of degree $2r$ and
of genus $r+1$ in $\bb P^{r}$.  So $C^\dg$ is of genus $r+1$.  Hence
  $\mc O_{C^\dg}(1)\simeq\omega$ by
Lemma~\ref{leextr}.  Hence the given isomorphism  $C\risom C^\dg$ is defined by
 $H^0(\omega)$, so coincides
with $\kappa$.  Thus the final assertion holds, and it implies (e).  The
proof is now complete.
\end{proof}

\begin{dfn}\label{dfBlwp}
For each integer $n\ge0$, set
 $$\omega^n := (\sSym^n\omega)\big/\!\Torsion(\sSym^n\omega).$$
 Plainly, $\bigoplus\omega^n$ is a quasi-coherent sheaf of finitely
 generated $\mc O_C$-domains.  Form
$$\wh C:= \Proj(\textstyle\bigoplus\omega^n)
\text{ and } \beta\:\wh C\to C,$$
 with $\beta$ the structure map.  In keeping with the notation $\mc O:=\mc
 O_C\text{ and } \?{\mc O}:= \nu_*\mc O_{\?C}$ and with the notation
 $\?\omega:=\nu_*\omega_{\?C}\text{ and } \?{\mc O}\omega := \nu_*(\mc
 O_{\?C}\omega)$, set
$$\wh{\mc O}:=\beta_*\mc O_{\wh C}\text{ and }
 \wh\omega:=\beta_*(\omega_{\wh  C})\text{ and }
\wh{\mc O}\omega := \beta_*(\mc O_{\wh  C}\omega).$$
 Call $\wh C$ the {\it blowup} of $C$ with respect to $\omega$, and
 $\beta$ the {\it blowup map}. 
\end{dfn}

\begin{prp}\label{prUP}
  The blowup $\wh C$ is an integral curve, the blowup map $\beta\:\wh
  C\to C$ is birational, and the sheaf $\mc O_{\wh C}\omega$ is
  invertible and generated by $H^0(\omega)$.  Furthermore, given any
  integral scheme $A$ and any nonconstant map $\alpha\:A\to C$, the
  sheaf $\mc O_A\omega$ is invertible iff there is a map $a\:A\to\wh
  C$ such that $\alpha=\beta a$; if so, then $\mc O_A\omega=a^*\mc
O_{\wh C}\omega$.
  \end{prp}
 \begin{proof}
  Plainly, $\wh C$ is an integral scheme, and $\beta\:\wh C\to C$ is of
  finite type.  Now, the smooth locus of $C$ is a nonempty open set $U$
  on which $\omega$ is invertible.  Fix any $U$ on which $\omega$ is
  invertible.  On $U$, the sum $\bigoplus\omega^n$ is locally isomorphic
  to the polynomial algebra in one variable over $\mc O_C$; whence,
  $\beta$ restricts to an isomorphism over $U$.  Thus $\beta\:\wh C\to
  C$ is birational, and so $\wh C$ is an integral curve.

  Since $\wh C$ is integral and $\mc O_{\wh C}(1)$ is invertible, $\mc
  O_{\wh C}(1)$ is torsion free.  Now, $\omega$ is invertible on $U$;
  so, on $\beta^{-1}U$, the tautological surjection $\beta^*\omega
  \onto\mc O_{\wh C}(1)$ is an isomorphism.  Hence, this surjection
  induces a global isomorphism $\mc O_{\wh C}\omega\risom \mc O_{\wh C}(1)$.
  Thus $\mc O_{\wh C}\omega$ is invertible.  And it is generated by
  $H^0(\omega)$ because $\omega$ is.

  Suppose $\alpha=\beta a$.  Then $\alpha^*\omega = a^*\beta^*\omega$.
  Now, the surjection $\beta^*\omega\onto\mc O_{\wh C}\omega$
  is an isomorphism on $\beta^{-1}U$.  Hence it induces a surjection
$\alpha^*\omega \onto
  a^*\mc O_{\wh C}\omega$, which is an isomorphism on $\alpha^{-1}U$; and
  $\alpha^{-1}U$ is nonempty since $\alpha$ is nonconstant.  Since $\mc
  O_{\wh C}\omega$ is invertible, so is $a^*\mc O_{\wh C}\omega$.  Hence
  $a^*\mc O_{\wh C}\omega$ is torsion free since $A$ is integral.
  Therefore, there is an induced isomorphism $\mc O_A\omega \risom
  a^*\mc O_{\wh C}\omega$.  Thus $\mc O_A\omega$ is invertible.

  Conversely, assume $\mc O_A\omega$ is invertible.  Then the surjection
  $\alpha^*\omega\onto\mc O_A\omega$ defines a $C$-map $p\:A\to\bb
  P(\omega)$  by general
  principles \cite[Prp.\,(4.2.3), p.\,73]{ega2}.  Plainly, $p$ factors
  through $\wh C\subset\bb P(\omega)$
  if, for each $n\ge0$, the induced map $w_n\:\alpha^*\sSym^n\omega\to(\mc
  O_A\omega)^{\ox n}$ factors through $\alpha^*\omega^n$.  But $w_n$ does
  factor because, on the one hand, the surjection
  $\alpha^*\sSym^n\omega\onto \alpha^*\omega^n$ is an isomorphism on
  $\alpha^{-1}U$, so its kernel is torsion, and on the other hand, $(\mc
  O_A\omega)^{\ox n}$ is torsion free since $\mc O_A\omega$ is
  invertible and $A$ is integral.
\end{proof}

\begin{dfn}\label{dfGrlcs}
  Denote by $G$ the largest open subset of $C$ where $\omega$ is
  invertible\emdash that is, where the local rings of $C$ are
  Gorenstein\emdash and call $G$ the {\it Gorenstein locus}.
\end{dfn}

\begin{cor}\label{coGrlcs}
The Gorenstein locus  $G$ of $C$ is the largest open subset  $A$ such
that the blowup map  $\beta$ restricts to an isomorphism
$\beta^{-1}A\risom A$. 
\end{cor}
\begin{proof}
  Taking $U:=G$ in the first paragraph of the proof of
  Proposition~\ref{prUP}, we find that $\beta$ restricts to an
  isomorphism $\beta^{-1}G\risom G$.  On the other hand, given an $A$
  such that $\beta$ restricts to an isomorphism $\beta^{-1}A\risom A$,
  plainly $\omega$ is invertible on $A$ as $\mc O_{\wh C}\omega$ is
  invertible on $\wh C$ by Proposition~\ref{prUP}; whence, $A$ lies in
  $G$.
\end{proof}

\begin{dfn}\label{dfCH} Owing to Proposition~\ref{prUP}, the
  normalization map $\nu\: \?C \to C$ and the canonical map $\?\kappa\:
  \?C\to C'$, both factor uniquely through normalization map $\wh\nu \:
  \?C \to \wh C$; that is,
$$\nu=\beta\wh\nu \text{ and }\?\kappa =\wh\kappa\wh\nu.$$
In view of Corollary~\ref{coGrlcs}, it is natural to set
$\kappa:=\wh\kappa\circ(\beta|G)^{-1}$, thereby extending the definition
of $\kappa$ for a Gorenstein $C$ to an arbitrary $C$.
In addition to $\?\kappa\: \?C\to C'$,  call
  $\wh\kappa\:\wh C\to C'$ and $\kappa\:G\to C'$ the {\it canonical
    maps\/} of $C$.
\end{dfn}

\begin{cor}\label{coDeg}
The invertible sheaf $\mc O_{\wh C}\omega$ is of degree $2g-2-\eta$. 
\end{cor}
\begin{proof}
  Proposition~\ref{prUP} implies $\wh\nu^*\mc O_{\wh C}\omega=\mc O_{\?
    C}\omega$; so $\mc O_{\wh C}\omega$ and $\mc O_{\?  C}\omega$ are of
  the same degree.  Hence Lemma~\ref{prDeg} yields the assertion.
\end{proof}


\begin{lem}\label{le3.1}
  Let $\mc F$ be a coherent sheaf on $C$.  Let $P\in C$ be a (closed)
  point at which $\mc F$ is invertible and generated by $H^0(\mc F)$.
  Assume either
\begin{enumerate}
\item there exists a (closed) point $Q\in C$ such that $Q\neq P$, such
  that $\mc F$ is invertible at $Q$ too, and such that
  $H^0(\sM_{\{Q\}}\sM_{\{P\}}\mc F)= H^0(\sM_{\{Q\}}\mc F)$, or
 \item  the natural map $v\:H^0(\sM_{\{P\}}\mc F)\to H^0(\sM_{\{P\}}\mc
  F/\sM_{\{P\}}^2\mc F)$ is not surjective.
\end{enumerate}
Then there exists a coherent subsheaf $\mc G\subset\mc F$ such that
	$$h^0(\mc G)=h^0(\mc F)-1\text{ and }h^1(\mc G)=h^1(\mc F)+1$$
and such that  $\Supp(\mc F/\mc G)$ is $\{P\}\cup \{Q\}$ if \(a) holds
or is  $\{P\}$ if \(b) holds.
\end{lem}
\begin{proof}
  Suppose (a) holds.  Take $\mc G:=\sM_{\{Q\}}\sM_{\{P\}}\mc F$.  Then
  $\Supp(\mc
  F/\mc G)$ consists of $P\cup Q$.  And $h^0(\mc F/\mc G)=2$ since $\mc
  F$ is invertible at both $P$ and $Q$; whence, $\chi(\mc G) = \chi(\mc F)
  -2$.  So it remains to prove $h^0(\mc G)=h^0(\mc F)-1$.

As $\mc F$ is invertible at $Q$, we have $h^0(\mc F/\sM_{\{Q\}}\mc
F)=1$.  So the inclusion $\sM_{\{Q\}}\mc F\into\mc F$ yields
$h^0(\sM_{\{Q\}}\mc F)\ge h^0(\mc F)-1$.  By hypothesis,
$H^0(\sM_{\{Q\}}\sM_{\{P\}}\mc F)= H^0(\sM_{\{Q\}}\mc F)$.  Hence
$h^0(\sM_{\{Q\}}\sM_{\{P\}}\mc F)\ge h^0(\mc F)-1$.  But
$H^0(\sM_{\{Q\}}\sM_{\{P\}}\mc F)\subset H^0(\sM_{\{P\}}\mc F)$.  And
$H^0(\sM_{\{P\}} \mc F) \subsetneqq H^0(\mc F)$ because $\mc F$ is
generated by $H^0(\mc F)$ at $P$.  Therefore,
$h^0(\sM_{\{P\}}\sM_{\{Q\}}\mc F)\le h^0(\mc F)-1$.  Thus $h^0(\mc
G)=h^0(\mc F)-1$, and the proof is complete when (a) holds.
 
Suppose (b) holds.  Take a vector subspace $V$ of $H^0(\sM_{\{P\}}\mc
F/\sM_{\{P\}}^2\mc F)$ such that $V$  contains the image of $H^0(\sM_{\{P\}}\mc
F)$ and  $V$  is of codimension 1.  Then take a subsheaf
$\mc H$ of $\sM_{\{P\}}\mc F/\sM_{\{P\}}^2\mc F$ such that $H^0(\mc
H)=V$.  Let $\mc G$ be the preimage of $\mc H$ in $\sM_{\{P\}}\mc F$.
Then $\Supp(\mc F/\mc G)=\{P\}$.  Now, $h^0(\mc F/\mc G)=h^0(\mc
F/\sM_{\{P\}}\mc F)+ h^0(\sM_{\{P\}}\mc F/\mc G)$.  But $h^0(\mc
F/\sM_{\{P\}}\mc F)=1$ as $\mc F$ is invertible at $P$, and
$h^0(\sM_{\{P\}}\mc F/\mc G)=1$ by construction.  So $h^0(\mc F/\mc
G)=2$; whence, $\chi(\mc G) = \chi(\mc F) -2$.  So it remains to prove
$h^0(\mc G)=h^0(\mc F)-1$.

Form this natural commutative diagram with exact rows:
$$\begin{CD}\label{CD1}
  0 @>>> H^0(\mc G)@>>> H^0(\sM_{\{P\}}\mc F)@>w>> H^0(\sM_{\{P\}}\mc F/\mc G)\\
  @.      @VVV            @VvVV  @|\\
  0 @>>> H^0(\mc G/\sM_{\{P\}}^2\mc F)@>u>> H^0(\sM_{\{P\}}\mc F/\sM_{\{P\}}^2\mc F)
  @>>> H^0(\sM_{\{P\}}\mc F/\mc G)
\end{CD}
$$
By construction, $\mc G/\sM_{\{P\}}^2\mc F=\mc H$.  Hence $\rIm(u)=V$.  By
construction, $V\supset\rIm(v)$.  Therefore, $w=0$.  Hence $H^0(\mc
G)\risom H^0(\sM_{\{P\}}\mc F)$.  However, $h^0(\mc F/\sM_{\{P\}}\mc
F)=1$.  And $H^0(\sM_{\{P\}} \mc F) \subsetneqq H^0(\mc F)$ because $\mc
F$ is generated by $H^0(\mc F)$ at $P$.  Therefore, the inclusion
$\sM_{\{P\}}\mc F\into\mc F$ gives $h^0(\sM_{\{P\}}\mc F)= h^0(\mc
F)-1$.  Thus $h^0(\mc G)=h^0(\mc F)-1$, and the proof is complete.
\end{proof}

\begin{lem}\label{le3hyp}
Let $\mc G\subset \omega$ be a coherent subsheaf.  Assume that
$\Supp(\omega/\mc G)$ lies in the Gorenstein locus $G$ and that
$h^0(\mc G)=g-1$ and $h^1(\mc G)=2$.  Then $C$ is hyperelliptic.
\end{lem}
\begin{proof}
Set $\mc L:=\sHom(\mc G,\,\omega)$.  Then $h^0(\mc L)=2$ by duality.  So
there is an $f$ in $H^0(\mc L)$ that is not a multiple of the inclusion
$h\:\mc G\into \omega$.  Set $\mc H:=\mc G+f\mc G\subset\omega$.  Then
$H^0(\mc G) \subseteq H^0(\mc H)$.

For a moment, suppose $H^0(\mc G) = H^0(\mc H)$.  Then $f$ induces a
$k$-linear endomorphism of $H^0(\mc G)$.  View $f$ as multiplication by
an element of the function field of $C$.  Then, by the Cayley--Hamilton
theorem, $f$ is integral over $k$, so lies in $k$, since $k$ is
algebraically closed.  Therefore, $f$ is a multiple of $h\:\mc G\into
\omega$, contrary to the choice of $f$.  Thus $H^0(\mc G) \subsetneqq
H^0(\mc H)$.

By hypothesis, $h^0(\mc G)=g-1$.  So $h^0(\mc H)=g$.  Hence $H^0(\mc H)=
H^0(\omega)$.  But $H^0(\omega)$ generates $\omega$.  Hence $\mc
H=\omega$.

Set $S:=\Supp(\omega/\mc G)$.  By hypothesis, $S\subset G$.  So $\omega$
is invertible along $S$.  It follows that $f\mc G=\omega$ along $S$,
since $\mc G+f\mc G=:\mc H$ and $\mc H=\omega$.  Hence $\mc G$ is
invertible along $S$.  But $\mc G$ is equal to $\omega$ off $S$.
Therefore, $\mc L$ is invertible.  Furthermore, $\mc L$ is generated by
its two global sections $f$ and $h$.

By duality, $h^0(\mc L)=2$ and $h^1(\mc L)=g-1$.  Hence, by the
Riemann--Roch Theorem, $\deg(\mc L)=2$.  Therefore, the pair $(\mc L,\,
H^0(\mc L))$ defines a map $C\to \bb P^1$ of degree 2.  Thus $C$ is
hyperelliptic, as asserted.
\end{proof}

\begin{lem}\label{leSep}
 \(1)   Let $P\in C$ be a (closed) point.  Then  $h^1(\sM_{\{P\}}\omega)=1$.

\(2) Let $P,\,Q\in C$ be distinct (closed) points, with $Q$ multiple.
Then
\begin{align*}
 &\text{\rm(a)}\enspace h^1(\sM_{\{P\}}\sM_{\{Q\}}\omega)=1, \text{ and}\\
 &\text{\rm(b)}\enspace h^0(\sM_{\{P\}}\sM_{\{Q\}}\omega) < 
\min\{\,h^0(\sM_{\{P\}}\omega),\,h^0(\sM_{\{Q\}}\omega)\,\}.
\end{align*}
\end{lem}
\begin{proof}
Consider (1).  First, suppose $P$ is simple.  Then $\omega$ is
invertible at $P$.  Hence the inclusion $\sM_{\{P\}}\omega\into\omega$
yields this long exact sequence: $$ H^0(\sM_{\{P\}}\omega)\xto{u}
H^0(\omega)\to k \to H^1(\sM_{\{P\}}\omega) \to H^1(\omega)\to 0.$$ Now,
$H^0(\omega)$ generates $\omega$. Hence $u$ isn't surjective.
So $H^1(\sM_{\{P\}}\omega) \risom H^1(\omega)$. But $h^1(\omega)=1$.
Thus (1) holds when $P$ is simple.

Suppose $P$ is multiple. Then $\sC\subset \sM_{\{P\}}$.  Hence
$\sC\omega\subset \sM_{\{P\}}\omega\subset \omega$.  But $\sC\omega
=\nu_*(\omega_{\?C})$ by Lemma~\ref{leCshf}; whence, $h^1(\sC\omega)=1$.
Therefore,
$$1=h^1(\sC\omega)\ge h^1(\sM_{\{P\}}\omega)\ge h^1(\omega)=1.$$
Thus (1) also holds when $P$ is multiple, and so (1) always holds.

Consider (2).  Again, first suppose $P$ is simple.  Then the
normalization map $\nu\: \?C \to C$ is an isomorphism over $P$.  Set
 $\?P:=\nu^{-1}P$, and let $\sM_{\{\?P\}}$ denote its maximal ideal sheaf.
Then $h^1(\sM_{\{\?P\}}\omega_{\?C})=1$ by (1) applied to $\?P\in\?C$.
However, $\nu_*\bigl(\sM_{\{\?P\}}\omega_{\?C}\bigr) =\sM_{\{P\}}\?\omega$,
and $\?\omega=\sC\omega$ by Lemma~\ref{leCshf}.  Hence
$h^1(\sM_{\{P\}}\sC\omega)=1$.

Since $Q$ is multiple, $\sC\subset \sM_{\{Q\}}$.  Hence
$\sM_{\{P\}}\sC\omega\subset \sM_{\{P\}}\sM_{\{Q\}}\omega\subset
\omega$.  Therefore,
$$1=h^1(\sM_{\{P\}}\sC\omega)\ge h^1(\sM_{\{P\}}\sM_{\{Q\}}\omega)\ge
h^1(\omega)=1.$$ Thus (2)(a) holds when $P$ is simple.

Suppose that $P$ is multiple.  Then $\sC\subset \sM_{\{P\}}\sM_{\{Q\}}$
since $P$ and $Q$ are distinct.  So $\sC\omega\subset
\sM_{\{P\}}\sM_{\{Q\}}\omega$.  But $\sC\omega =\nu_*(\omega_{\?C})$ by
Lemma~\ref{leCshf}, and $\sM_{\{P\}}\sM_{\{Q\}}\omega \subset \omega$.
Therefore, $$1=h^1(\sC\omega)\ge h^1(\sM_{\{P\}}\sM_{\{Q\}}\omega)\ge
h^1(\omega)=1.$$ Thus (2)(a) also holds when $P$ is multiple, and so
(2)(a) always holds.

Finally, the inclusion $\sM_{\{P\}}\sM_{\{Q\}}\omega\into
\sM_{\{P\}}\omega$ yields this long exact sequence:
\begin{align*}\label{al3.8.1}
0&\to H^0(\sM_{\{P\}}\sM_{\{Q\}}\omega)\to H^0(\sM_{\{P\}}\omega)
\to H^0(\sM_{\{P\}}\omega/\sM_{\{P\}}\sM_{\{Q\}}\omega)\\
&\to H^1(\sM_{\{P\}}\sM_{\{Q\}}\omega)\to H^1(\sM_{\{P\}}\omega)\to0.
\end{align*}
Hence (1) and (2)(a) imply $h^0(\sM_{\{P\}}\sM_{\{Q\}}\omega)<
h^0(\sM_{\{P\}}\omega)$.  Similarly, the inclusion
$\sM_{\{P\}}\sM_{\{Q\}}\omega\into \sM_{\{P\}}\omega$ yields
$h^0(\sM_{\{P\}}\sM_{\{Q\}}\omega)< h^0(\sM_{\{Q\}}\omega)$ as
$h^1(\sM_{\{Q\}}\omega)=1$ by (1) with $P:=Q$.  Thus (2)(b) holds, and
the proof is complete.
\end{proof}

\begin{thm}\label{prGiso}
  If $C$ is nonhyperelliptic, then $\kappa\:G\to C'$ is an open
  embedding.
\end{thm}
\begin{proof}
Set $\wh G:=\beta^{-1}G$.  Then $\beta$ restricts to an isomorphism
$\beta^{-1}G\risom G$ by Corollary~\ref{coGrlcs}.  So it suffices to
prove $\wh\kappa$ restricts to an isomorphism $\wh G\risom
\wh\kappa\wh G$.

 Let $\wh P\in \wh G$ be an arbitrary (closed) point.  For a moment,
assume (a) that $H^0(\omega)$ separates $\wh P$ from every other point
$\wh Q\in \wh C$ and (b) that $H^0(\omega)$ separates tangent vectors at
$\wh P$.  More precisely, (a) means that there is an $f\in H^0(\omega)$
whose image in $H^0(\mc O_{\wh C}\omega)$ lies in $H^0(\sM_{\{\wh Q\}}\mc
O_{\wh C}\omega)$, but not in $H^0(\sM_{\{\wh P\}}\mc O_{\wh C}\omega)$.
And (b) means that $H^0(\sM_{\{\wh P\}}\mc O_{\wh C}\omega)$ maps onto
$H^0(\sM_{\{\wh P\}}\mc O_{\wh C}\omega\big/\sM_{\{\wh P\}}^2\mc O_{\wh
C}\omega)$.

Condition (a) implies that $\wh\kappa\wh Q\not=\wh\kappa\wh P$, because
$f$ defines a hyperplane in $\bb P^{g-1}$ that contains $\wh\kappa\wh
Q$, but not $\wh\kappa\wh P$.  Taking $\wh Q\in \wh G$ shows that the
restriction $\wh G\to C'$ is injective.  Taking $\wh Q\notin \wh G$
shows that $\wh G= \wh\kappa^{-1}\wh\kappa\wh G$; whence, $\wh G\to\wh
\kappa\wh G$ is finite since $\wh\kappa\:\wh C\to C'$ is finite.
Condition (b) now implies that $\wh G\to \wh\kappa\wh G$ is an isomorphism
owing to a simple lemma of Commutative Algebra \cite[Lem.\,II-7.4,
p.\,153]{H}.  Thus it remains to prove (a) and (b).

Set $P:=\beta\wh P$ and $Q:=\beta\wh Q$.  To prove (a), it suffices to
prove the inequality
\begin{equation}\label{eqprGiso1}
  h^0(\sM_{\{P\}}\sM_{\{Q\}}\omega)< h^0(\sM_{\{Q\}}\omega).
\end{equation}
Indeed, (\ref{eqprGiso1}) implies that there is an $f $ in
$H^0(\sM_{\{Q\}}\omega)$ not in $H^0(\sM_{\{P\}}\omega)$.  But
$\sM_{\{Q\}}\mc O_{\wh C}\subset \sM_{\{\wh Q\}}$, and $\beta\:\wh C\to C$
is an isomorphism at $\wh P$ by Corollary~\ref{coGrlcs}.  Hence (a)
holds.  Now, (\ref{eqprGiso1}) holds if $Q$ is multiple by
Lemma~\ref{leSep}(2), and if $Q$ is simple by Lemma~\ref{le3.1} with
$\mc F:=\omega$ and by Lemma~\ref{le3hyp}.

Finally, $H^0(\sM_{\{P\}}\omega)$ maps onto $ H^0(\sM_{\{P\}}\omega\big/
\sM_{\{P\}}^2\omega)$ owing to Lemma~\ref{le3.1} with $\mc F:=\omega$ and to
Lemma~\ref{le3hyp}. But $\beta\:\wh C\to C$ is an isomorphism at $\wh P$
by Proposition~\ref{coGrlcs}.  Thus (b) holds, and the proof is
complete.
\end{proof}

\section{Arithmetically normal models}\label{sc_nl}
 Preserve the general setup introduced at the beginning of
 Section~\ref{sc_cm}, after Proposition~\ref{prHyp2}, and before
 Definition~\ref{dfBlwp} and Definition~\ref{dfCH}.  In this section,
 the main result is Theorem~\ref{thAn}, which gives necessary and
 sufficient conditions for the canonical model $C'$ to be arithmetically
 normal.  The main new tool is Castelnuovo Theory.   It also yields a
 third proof, given in Remark~\ref{rm3rdpf}, that, if $C$ is
 nonhyperelliptic and Gorenstein, then the
 canonical map yields an isomorphism $\kappa\:C\risom C'$.

\begin{lem}\label{leAmp}
  Let $\mc F$ be a torsion-free sheaf of rank $1$ on $C$, and $P\in C$.

\(1)  If $\chi(\mc F)\ge g$, then $h^1(\mc F)=0$. 

\(2) If $\chi(\mc F)\ge g+h^0(\mc F/\sM_{\{P\}}\mc F)$, then $H^0(\mc F)$
generates $\mc F$ at  $P$.

\(3) If $\mc F$ is invertible and $\deg(\mc F)\ge 2g+1$, then $\mc F$ is
very ample and $h^1(\mc F)=0$.
 \end{lem}
 \begin{proof}
To prove (1), suppose $h^1(\mc F)\not=0$.  Then, by duality, there is a
nonzero map $\mc F\to \omega$.  It is injective as $\mc F$ is
torsion-free of rank $1$.  Hence $\chi(\mc F)\le\chi(\omega)$, contrary
to hypothesis.  Thus (1) holds.

To prove (2), note that $\chi(\sM_{\{P\}}\mc F)=\chi(\mc F)-h^0(\mc F/
\sM_{\{P\}}\mc F)\ge g$; so (1) implies that $h^1(\sM_{\{P\}}\mc F)=0$.  Hence
$H^0(\mc F)\to H^0(\mc F/\sM_{\{P\}}\mc F)$ is surjective.  Thus (2) holds.

To prove (3), note that $\chi(\mc F)\ge g+2$ by the Riemann--Roch
theorem.  Hence (1) implies that $h^1(\mc F)=0$, and (2) implies that
$H^0(\mc F)$ generates $\mc F$ at every point.  Suppose $\mc F$ is not
very ample.  Then either $H^0(\mc F)$ does not separate points or it
does not separate tangent directions; either way, Lemma~\ref{le3.1}
implies that there is a coherent subsheaf $\mc G\subset\mc F$ such that
$h^0(\mc G)=h^0(\mc F)-1$ and $h^1(\mc G)=h^1(\mc F)+1$.  Then $h^1(\mc
G)\ge1$, and $\chi(\mc G)=\chi(\mc F)-2\ge g$, contrary to (1) applied
with $\mc F:=\mc G$.  Thus (3) holds, and the proof is complete.
  \end{proof}

\begin{prp}\label{prVA} Assume $C$ is not Gorenstein.  Then  $\mc O_{\wh
    C}\omega$ is very ample, and $h^1\bigl(\mc O_{\wh C}\omega^{\ox
    n}\bigr)=0$ for $n\ge1$.
 \end{prp}
\begin{proof}
Fix a (closed) point $P\in C$.  Set $\xi_P := \dim(\wh{\mc
O}_P/\mc O_P)$.  Recall $\wh{\mc O}\omega := \beta_*(\mc O_{\wh
C}\omega)$.  Set $\mu_P:=\dim((\wh{\mc O}\omega)_P/\omega_P)$.
Lemma~\ref{leCshf} yields an $x\in\omega_P$ so that $\?{\mc O}_Px =
(\?{\mc O}\omega)_P$.  Then $\wh{\mc O}_Px = (\wh{\mc O}\omega)_P$ since
$\mc O_{\wh C}\omega$ is invertible.  Hence $\dim((\wh{\mc O}\omega)_P/\mc
O_Px)=\xi_P$.  Hence
 $$\xi_P=\eta_P+\mu_P$$ owing to Equation~\ref{eqEta1}.  If $C$ is
not Gorenstein at $P$, then $\eta_P\ge1$ by Lemma~\ref{leEta}, and
$\mu_P\ge1$ by  the proof of Proposition 28 in \cite[p.\,438]{BF}.

Let $\wh g$ denote the arithmetic genus of $\wh C$.  Then $\wh g=g-
\sum_{P\in C}\xi_P$.  Set $\mu:=\sum\mu_P$.  Then $\wh g=g-
\eta-\mu$.  However, $\deg\mc O_{\wh C}\omega= 2g-2-\eta$ by
Corollary~\ref{coDeg}.  Therefore, $\deg\mc O_{\wh C}\omega= 2\wh
g-2+\eta+2\mu$.  But $\eta\ge1$ and $\mu\ge1$ since $C$ is not
Gorenstein.  Hence $\deg\bigl(\mc O_{\wh C}\omega^{\ox n}\bigr)\ge 2\wh
g+1$ for $n\ge1$.  So Lemma~\ref{leAmp}(3) yields the assertion.
  \end{proof}

\begin{dfn}\label{dfLnPn}
  Call the canonical model $C'\subset\bb P^{g-1}$ {\it linearly
    normal\/} if the linear series of hyperplane sections is complete, in other
  words, if $h^0(\mc O_{C'}(1))=g$.

  Call $C'$ {\it projectively normal\/} if the linear series of
  hypersurface sections of degree $n$ is complete for every $n\ge1$, in
  other words, if  the natural map
$$\Sym^nH^0(\omega)\to H^0(\mc O_{C'}(n))$$
is surjective for every $n\ge1$.
\end{dfn}

\begin{lem}\label{leLnPn}
\(1) The model  $C'$ is linearly normal iff it is projectively normal.

\(2) If $C$ is Gorenstein, then $C'$ is linearly normal.

\(3) If $C$ is not Gorenstein, then $h^1(\mc O_{C'}(l))=0$ for $l\ge1$;
furthermore, then $C'$ is linearly normal iff $d'=g+g'-1$.
 \end{lem}
 \begin{proof}
Trivially, $C'$ is linearly normal if it is projectively normal.  In
other words, sufficiency holds in (1).

If $C$ is Gorenstein,  then $\kappa\:C\to C'$ exists and
$\kappa^*\mc O_{C'}(1)=\omega$.  Hence, then $h^0(\mc O_{C'}(1))\le
h^0(\omega)=g$.  But $h^0(\mc O_{C'}(1))\ge g$ always.  Thus (2) holds.

Assume $C$ is hyperelliptic.  Then $C'$ is equal to the rational
normal curve of degree $g-1$ in $\bb P^{g-1}$ by Theorem~\ref{coRnc}
or by Proposition~\ref{prHyp2}(3).  So $g'=0$ and $h^1(\mc
O_{C'}(1))=0$; furthermore, $C'$ is linearly normal and projectively
normal.  Thus (1) and (2) hold in this case.  And (3) does not apply,
since $C$ is Gorenstein by Proposition~\ref{prHyp2}(2).  So from now
on, assume $C$ is nonhyperelliptic.

To prove necessity in (1) and to prove (3),
let's apply Castelnuovo Theory as presented in
\cite[pp.\,114--117]{ACGH}.  First, note that, by the General Position
Theorem \cite[p.\,109]{ACGH}, a general hyperplane $H$ meets $C'$ in a
set $\Gamma$ of $d'$ distinct points any $g-1$ of which are linearly
independent.  Second, for $l\ge1$, form the linear series of
hypersurfaces of degree $l$ in $H$ containing $\Gamma$, and note that,
by the lemma on p.\,115 in \cite{ACGH}, the series has (projective)
dimension at least $\min\{d'-1,\ l(g-2)\}$.

Consider the following standard left exact sequence:
\begin{equation}\label{eqLP3}
 0\to H^0(\mc O_{C'}(l-1))\xto{u} H^0(\mc O_{C'}(l))\xto{v} H^0(\mc
O_{\Gamma}(l)). 
 \end{equation}
 Let $V_l$ denote the image of $H^0(\mc O_{\bb P^{g-1}}(l))$ in $H^0(\mc
O_{C'}(l))$, and set $W_l:=v(V_l)$.  Then $\dim(W_l)\ge \min\{d',\ l(g-2)+1\}$
by the second note above.  Hence
 \begin{equation}\label{eqLP1}
 h^0(\mc O_{C'}(l))-h^0(\mc O_{C'}(l-1))\ge\min\{d',\ l(g-2)+1\}.
\end{equation}

Also, if equality holds in (\ref{eqLP1}), then $v\bigl(H^0(\mc
O_{C'}(l))\bigr)=W_l$ since both sides have the same dimension.  So
$H^0(\mc O_{C'}(l))$ is spanned by $V_l$ and ${\rm Im}(u)$.  But
$u(V_{l-1})\subset V_l$.  And, if $C'$ is linearly normal, then
$V_1=H^0(\mc O_{C'}(1))$.  Hence, if in addition, equality holds in
(\ref{eqLP1}) for $l\ge2$, then induction on $l$ yields $V_l=H^0(\mc
O_{C'}(l))$ for $l\ge1$; in other words, then $C'$ is projectively
normal.  Thus to complete the proof of (1), we have to prove that
equality holds in (\ref{eqLP1}) for $l\ge2$.

Set $h(l):=h^1(\mc O_{C'}(l))$.  Then $h^0(\mc O_{C'}(l))=ld'+1-g'+h(l)$
by the Riemann--Roch Theorem.  Hence the bound~(\ref{eqLP1}) is
equivalent to this bound:
\begin{equation}\label{eqLP2}
  d'-(h(l-1)-h(l))\ge\min\{d',\ l(g-2)+1\}.
\end{equation}
Here $h(l-1)-h(l)\ge0$ because the sequence~(\ref{eqLP3}) continues,
ending with $$H^1(\mc O_{C'}(l-1))\to H^1(\mc O_{C'}(l))\to 0.$$

Since $C$ is nonhyperelliptic, $d'=2g-2-\eta$ 
by Proposition~\ref{prdegC'}.  So
\begin{equation*}\label{eqeqLP4}
(l(g-2)+1)-d' = (l-2)(g-2)-1+\eta.
\end{equation*}
 For $l\ge2$, the right side is nonnegative unless $\eta=0$ and either
$l=2$ or $g=2$.  But Lemma~\ref{leEta} implies $\eta=0$ iff $C$ is
Gorenstein.  So, since $C$ is nonhyperelliptic, by
Proposition~\ref{prg=2}, the right side of (\ref{eqLP2}) is nonnegative
unless $\eta=0$ and $l=2$.

Hence, for $l\ge3$, the right side of (\ref{eqLP2}) is equal to $d'$.
But $h(l-1)-h(l)\ge0$.  Therefore, equality holds in (\ref{eqLP2}), and
$h(l-1)=h(l)$.  But, by Lemma~\ref{leAmp}(1) or by Serre's Theorem,
$h(l)=0$ for $l\gg0$.  So $h(l)=0$ for $l\ge2$.

Suppose $\eta>0$.  Then similarly, equality holds in (\ref{eqLP2}) for
$l=2$ too, and $h(1)=0$.  So equality holds in (\ref{eqLP1}) for
$l\ge2$, as desired.  Thus (1) holds.  And $h^0(\mc O_{C'}(1))=d'+1-g'$
by the Riemann--Roch Theorem.  Thus (3) holds.

Finally, instead suppose $\eta=0$.  Then $d'=2g-2$.  And $h^0(\mc
O_{C'}(1)) = g$ by (2).  So the Riemann--Roch Theorem yields $g=(2g-2)
+1-g'+h(1)$.  But $g'\ge g$ as $\kappa\:C\to C'$ exists.  Hence
$h(1)\ge1$.  Now, take $l=2$ in (\ref{eqLP2}), getting $2g-2-h(1)$ on
the left as $h(2)=0$, and $2g-3$ on the right.  Hence, $h(1)=1$ and
equality holds in (\ref{eqLP2}) for $l=2$.  So equality holds in
(\ref{eqLP1}) for $l\ge2$, as desired.  Thus (1) always holds.  The
proof is now complete.
  \end{proof}

\begin{prp}\label{prGpn}
If $C$ is Gorenstein, then $C'$ is projectively normal.
\end{prp}
\begin{proof}
The assertion is immediate from Lemma~\ref{leLnPn}(1),\,(2). 
\end{proof}

\begin{rmk}\label{rm3rdpf}
  Assume $C$ is nonhyperelliptic and Gorenstein.  Then Castelnuovo
  Theory yields a third proof that $\kappa\:C\to C'$ is an isomorphism
  (the other two are the proof of Theorem~\ref{thNonh} by computing
  global invariants and the proof of Theorem~\ref{prGiso}  by
  separating points and tangent directions.

Indeed, take a nonzero section in $ H^0(\omega)$, and use it to form,
for each $l\ge2$, the rows in the following commutative diagram:
$$\begin{CD}\label{CD2}
  0 @>>> H^0(\omega^{\ox(l-1)}) @>>> H^0(\omega^{\ox l})\\
  @.      @Aw_{l-1}AA           @Aw_lAA  \\
  0 @>>> H^0(\mc O_{C'}(l-1)) @>>> H^0(\mc O_{C'}(l))
\end{CD}
$$
The $w$'s are induced by $\kappa$, so are injective.  We have to prove
they are bijective, as $C=\Proj\bigl(\bigoplus H^0(\omega^{\ox l})\bigr)$ since
$\omega$ is ample, and as $C'=\Proj\bigl(\bigoplus H^0(\mc O_{C'}(l))\bigr)$.

Let's proceed by induction on $l$.  First off, $w_1$ is bijective
because $C'$ is linearly normal by Lemma~\ref{leLnPn}(2).  Suppose
$w_{l-1}$ is bijective.  Then $w_l$ is bijective iff
\begin{equation}\label{eq3.16a}
h^0(\omega^{\ox l}) - h^0(\omega^{\ox(l-1)}) =
h^0(\mc O_{C'}(l)) - h^0(\mc O_{C'}(l-1)).
\end{equation}
Now, Proposition~\ref{prBir} implies $\deg\kappa=1$; whence,
$\deg(\omega^{\ox l})=\deg(\mc O_{C'}(l))$.  Hence, owing to the Riemann--Roch
Theorem, Equation~\ref{eq3.16a} holds iff
\begin{equation}\label{eq3.16b}
h^1(\omega^{\ox l}) - h^1(\omega^{\ox(l-1)}) =
h^1(\mc O_{C'}(l)) - h^1(\mc O_{C'}(l-1)).
\end{equation}
But $h^1(\omega^{\ox l})=0$ for $l\ge2$ by Lemma~\ref{leAmp}, and
$h^1(\omega)=1$ by duality.  Furthermore, it was shown in the course of
the proof of Lemma~\ref{leLnPn} that $h^1(\mc O_{C'}(l))=0$ for $l\ge2$
and that $h^1(\mc O_{C'}(1))=1$.  Thus Equation~\ref{eq3.16b} holds, and
the proof is complete.
\end{rmk}

\begin{dfn}\label{dfnlyG}
  Call $C$ {\it nearly Gorenstein\/} if the complement of the Gorenstein
  locus consists of a single point $P$ and if the local ring $\mc O_P$
  is almost Gorenstein in the sense of Barucci and Fr\"oberg
  \cite[p.\,418]{BF}, namely, if
 $$\eta_P:=\dim(\?{\mc O}_P/\mc O_P) - \dim (\mc O_P/\sC_P) = \dim(\Ext^1(k,
 \mc O_P))-1$$
where $k$ is the algebraically closed ground field.
\end{dfn}

\begin{lem}\label{prNlyG1}
Assume $C$ is not Gorenstein.  Then these  conditions are equivalent:
\begin{enumerate}
 \item  $h^0(\mc O_{\wh C}\omega)=g$;
 \item  $H^0(\omega)=H^0(\mc O_{\wh C}\omega)$;
 \item  $C$ is nearly Gorenstein.
 \end{enumerate}
 If \(a)--\(c) hold, then $\wh\kappa\:\wh C\to C'$ is an isomorphism, and $C'$
is linearly normal.
 \end{lem}

 \begin{proof}
Note that $h^0(\omega)=g$ and $H^0(\omega)\subset H^0(\mc O_{\wh
C}\omega)$.  So (a) and (b) are equivalent.

 Preserve the notation of the proof of Proposition~\ref{prVA}.  Form the
exact sequence $$0\to\omega\to\wh{\mc O}\omega\to\wh{\mc
O}\omega/\omega\to0.$$ It yields $g-1+\mu=h^0(\mc O_{\wh C}\omega)$
since $h^1(\mc O_{\wh C}\omega)=0$ by that proposition.  Hence $h^0(\mc
O_{\wh C}\omega)=g$ iff $\mu=1$.  But $\mu:=\sum\mu_P$ and $\mu_P\ge0$.
Hence $\mu=1$ iff there is one and only one $P$ such that $\mu_P=1$.
But $\mu_P=1$ iff $\mc O_P$ is almost Gorenstein, but not Gorenstein, by
Proposition 28 in \cite[p.\,438]{BF}.  Thus (a) holds iff (c) holds.

The sheaf $\mc O_{\wh C}\omega$ is very ample by Proposition~\ref{prVA};
so $H^0(\mc O_{\wh C}\omega)$ defines an embedding of $\wh C$ into
projective space.  Assume also that $H^0(\omega)=H^0(\mc O_{\wh
  C}\omega)$.  Then this embedding is essentially the canonical map $\wh
\kappa\:\wh C\to C'$.

Alternatively, we can prove that $\wh\kappa$ is an isomorphism via
Castelnuovo Theory proceeding as in Remark~\ref{rm3rdpf}, but with $\wh
C$ and $\mc O_{\wh C}\omega$ in place of $C$ and $\omega$.  This time,
$h^1\bigl(\mc O_{\wh C}\omega^{\ox l}\bigr)=0$ for $l\ge1$ by
Proposition~\ref{prVA}; furthermore, it was shown in the course of the
proof of Lemma~\ref{leLnPn} that now $h^1(\mc O_{C'}(l))=0$ for $l\ge1$.

Since $H^0(\omega)=H^0(\mc O_{\wh C}\omega)$ and since $\wh\kappa$ is an
isomorphism, $h^0(\mc O_{C'}(1))=g$; in other words, $C'$ is linearly
normal, and the proof is complete.
 \end{proof}

\begin{dfn}\label{dfAn}
  Say that $C'$ is {\it arithmetically normal\/} if its homogeneous
  coordinate ring is normal.
\end{dfn}

\begin{thm}\label{thAn} If $C$ is not Gorenstein, then these seven conditions
are equivalent:
 \begin{enumerate}
 \item  $C'$ is arithmetically normal;
 \item  $C'$ is smooth and projectively normal;
 \item  $C'$ is smooth and linearly normal;
 \item  $C'$ is smooth and extremal;
 \item  $d'=g+\?g-1$;
 \item $C$  is nearly normal;
 \item  $C$ is nearly Gorenstein, and  $\smash{\wh C}$ is smooth.
 \end{enumerate}
If these conditions hold, then at its unique multiple point, $C$ is of
multiplicity $g-\?g+1$ and of (maximal) embedding dimension $g-\?g+1$.
  \end{thm}
\begin{proof}
 Let $A$ be the homogeneous coordinate ring of $C'$.  Then $A$ is normal
iff $C'$ is smooth and $A$ is of depth 2 at the irrelevant ideal by
Serre's criterion.  But
this depth condition holds iff $C'$ is projectively normal by a
well-known theorem
due to Grothendieck  \cite[(2.2.4)]{KL71}.  Thus (a)
and (b) are equivalent.

Conditions (b) and (c) are equivalent by Lemma~\ref{leLnPn}.

Suppose (c) holds.  Then $C'$ is linearly normal and $C$ is not
Gorenstein; hence, Lemma~\ref{leLnPn}(3) yields $d'=g+g'-1$.  In
addition, $C'$ is smooth.  Hence $\?\kappa\:\?C\to C'$ is an
isomorphism.  So $g'=\?g$.  Thus (e) holds.

Conversely, suppose (e) holds.  Then the Riemann--Roch Theorem yields
$$h^0(\mc O_{C'}(1))=(g+\?g-1)+1-g'+h^1(\mc O_{C'}(1)).$$
 Since $C$ is not Gorenstein, $h^1(\mc O_{C'}(1))=0$ by
Lemma~\ref{leLnPn}(3).  And $g\le h^0(\mc O_{C'}(1))$ as $C'\subset\bb
P^{g-1}$ is nondegenerate.  Therefore, $g\le g+\?g-g'$, so $g'\le \?g$;
whence, $\?\kappa\:\?C\to C'$ is an isomorphism.  So $C'$ is smooth, and
$g'=\?g$.  The latter yields $d'=g+g'-1$; so $C'$ is linearly
normal by Lemma~\ref{leLnPn}(3).  Thus (c) holds.

Conditions (e) and (f) are equivalent by Proposition~\ref{prLbd}.  Thus
(a)--(c) and (e)--(f) are equivalent.

Again, suppose (c) holds.  Then $C'$ is smooth.  Hence $\wh\kappa\:\wh
C\to C'$ is an isomorphism.  So $\wh C$ is smooth.  In addition, $C'$ is
linearly normal, or $h^0(\mc O_{C'}(1))=g$.  But $\wh\kappa^*\mc
O_{C'}(1)=\mc O_{\wh C}\omega$.  So $h^0(\mc O_{\wh C}\omega)=g$.  By
hypothesis, $C$ is not Gorenstein.  By Lemma~\ref{prNlyG1}, therefore,
$C$ is nearly Gorenstein.  Thus (g) holds.

Conversely, suppose (g) holds.  Then $C$ is nearly Gorenstein.  By
Lemma~\ref{prNlyG1}, therefore, $\wh\kappa\:\wh C\to C'$ is an
isomorphism, and $C'$ is linearly normal.  In addition, $\wh C$ is
smooth.  So $C'$ is smooth.  Thus (c) holds.  Thus (a)--(c),\,(e)--(g)
are equivalent.

The last assertion concerns the unique multiple point $P\in C$.  Denote
its maximal ideal in its local ring $\mc O_P$ by $\mf m$.  Then by 
general principles, the multiplicity of $P$ is just $h^0\bigl(\?{\mc O}_P
/\mf m\?{\mc O}_P\bigr)$; moreover, $h^0\bigl(\?{\mc O}_P /\mc
O_P\bigr)=g-\?g$.  But $\sM_{\{P\}}=\sC$ by (f); so $\mf m\?{\mc
  O}_P=\mf m$.  Thus the multiplicity is $g-\?g+1$ at $P$.

By (f) and (g), the ring $\mc O_P$ is almost Gorenstein, but not
Gorenstein.  Hence the endomorphism ring of $\smash{\sM_{\{P\}}}$ is
equal to $\wh{\mc O}_P$ by Proposition 28 in \cite[p.\,438]{BF}.  But
$\wh{\mc
  O}_P$ is smooth, so Gorenstein.  Hence $\mc O_P$ is of maximal
embedding dimension by Proposition 25 in \cite[p.\,436]{BF}.

It remains to prove that (d) is equivalent to the other conditions.
Clearly, we may assume $C'$ is smooth.  Then $\?\kappa\:\?C\to C'$ is an
isomorphism, and so $g'=\?g$.  Hence (e) holds if and only if $g'=d'-r$
with $r:=g-1$.

Since $C$ is not Gorenstein, Lemma~\ref{leDeg} yields $d'<2r$.  Hence,
by Lemma~\ref{leextr}, any curve of degree $d'$ in $\bb P^r$ is of genus
at most $d'-r$.  But there exist curves in $\bb P^r$ of degree $d'$ and
genus exactly $d'-r$; see Example~\ref{rmCstr} below.  Hence any
extremal curve of degree $d'$ in $\bb P^r$ is of genus $d'-r$.
Therefore, $C'$ is extremal if and only if (e) holds.  The proof is now
complete.
 \end{proof}

\begin{eg}\label{rmCstr}
  In Theorem~\ref{thAn}, if (a)--(g) hold, then 
 the preimage of the unique multiple point is an 
  effective divisor on $\?C$ of degree $g-\?g+1$.  Conversely, Serre
  \cite[Chap.\,IV, n$^{\text o} 4$, p.\,70]{Sr} explains how to construct
  such a  $C$: start with any smooth curve $\?C$ of genus
  $\?g$ and with any effective divisor $D$ on $\?C$ of degree $g-\?g+1$
  and then contract $D$.
\end{eg}

\section{Projectively normal models}\label{sc_Rt}
Preserve the general setup introduced at the beginning of
 Section~\ref{sc_cm}, after Proposition~\ref{prHyp2}, and before
 Definition~\ref{dfBlwp} and Definition~\ref{dfCH}.  In this section, we
 prove Rosenlicht's Main Theorem, \cite[Thm.\,17, p.\,189]{R}, which
 essentially asserts that, if $C$ is nonhyperelliptic, then the
 canonical map $\wh\kappa\:\wh C\risom C'$ is an isomorphism.  We then
 apply this result to characterize the non-Gorenstein curves $C$ whose
 canonical model $C'$ is projectively normal.  We begin by proving two
 lemmas.

\begin{lem}\label{leR1} Given an $x\in H^0(\omega)$, set $W_x:=\{\,f\in
  k(C)\mid fx\in H^0(\omega)\,\}$.  Then, given a point $P\in C$, there
  exists an $x\in H^0(\omega)$ with these four properties:
\begin{enumerate}
 \item $(\wh{\mc O}\omega)_P = \wh{\mc O}_Px$;
 \item $\wh{\mc O}_P=\mc O_P[W_x]$;
 \item $k[W_x]$ is the ring of an affine open subset of $C'$, which
contains $\wh\kappa\beta^{-1}P$; 
\item $\mc O_P \subset W_x + \sC_P$.
\end{enumerate} 
\end{lem}
\begin{proof}
 By Proposition~\ref{prUP}, the sheaf $\mc O_{\wh C}\omega$ is
invertible and generated by $H^0(\omega)$.  So, as the base field $k$ 
is infinite, there exists an $x\in H^0(\omega)$ with Property~(a).

Let's now prove that any $x\in H^0(\omega)$ having (a) also has
(b)--(d).  First, let's prove that any $y\in\omega_P$ can be expressed
as a sum $y=y'+y''$ with $y'\in H^0(\omega)$ and $y''\in \sC_Px$.
Indeed, form the long exact sequence
 \begin{equation}\label{eqR1a}
H^0(\omega) \xto{u} H^0(\omega/\?\omega)\to H^1(\?\omega)
 \to H^1(\omega)\to 0.
  \end{equation}
  Plainly, $h^1(\omega)=1$ and $h^1(\?\omega)=1$.  Hence $u$ is
  surjective.  But it is clear that $H^0(\omega/\?\omega)=
  \bigoplus_Q(\omega_Q/\?\omega_Q)$.  So there is a $y''\in \?\omega_P$
  such that $y-y''\in H^0(\omega)$.  But $\?\omega_P=\sC_Px$ by
  Lemma~\ref{leCshf} as $\?{\mc O}_Px= (\?{\mc O}\omega)_P$ owing to (a).  So
  $y''\in \sC_Px$.

To prove (b), set $$V_x:=\{\,f\in k(C)\mid fx\in \omega_P\,\}.$$  Then
$\wh{\mc O}_P= \mc O_P[V_x]$ owing to the construction of $\wh C$ and to
(a).  Now, given $f\in V_x$, take $y:=fx$ above.  Say $y=y'+y''$ with
$y'\in H^0(\omega)$ and $y''\in \sC_Px$.  Then $y'=f'x$ with $f'\in
W_x$, and $y''=f''x$ with $f''\in\sC_P\subset \mc O_P$.  And $f=f'+f''$.
Therefore, $\mc O_P[V_x]=\mc O_P[W_x]$.  Thus (b) holds.

To prove (c), note that, by construction, $C'=\bigcup_{v\in H^0(\omega)}
\Spec\bigl(k[W_v]\bigr)$.  Now, $x$ vanishes nowhere on $\beta^{-1}P$
owing to (a).  Hence (c) holds.

Finally, to prove (d), let $f\in\mc O_P$.  Then $f\in V_x$.  So the
proof of (b) yields a decomposition $f=f'+f''$ with $f'\in W_x$, and
$f''\in\sC_P$.  Thus (d) holds.
  \end{proof}

\begin{rmk}\label{rmSt}
 St\"ohr \cite[Thm.\,3.2, p.\,123]{S} proved that $H^0(\omega)$
generates $\omega$ by introducing the ideas used in the proof of
Lemma~\ref{leR1} and developing them essentially as follows.  Consider
the map $u$ in the sequence~(\ref{eqR1a}).  It is surjective, and its
target is equal to $\bigoplus_Q(\omega_Q/{\sC}_Q\omega_Q)$ as
$\?\omega=\sC\omega$ by Lemma~\ref{leCshf}.  Hence, by Nakayama's
lemma, $H^0(\omega)$ generates $\omega_Q$ when $Q$ is multiple.

Finally, when $Q$ is simple, then $H^0(\omega)$ generates $\omega_Q$ by
the usual argument.  Namely, $h^0(\mc O_C(Q))=1$ since $g>0$.  So
$h^1(\omega(-Q))=1$ by duality.  Hence $H^0(\omega)\to
H^0(\omega/\omega(-Q))$ is surjective, as desired.
  \end{rmk}

\begin{lem}\label{leR2}
  Let $\mc F$ be a nonzero proper ideal of $\?{\mc O}$.  View
  $H^0(\sHom(\mc F,\,\omega))$ and $H^0(\omega)$ as subsets of
  $H^0(\sHom(\mc F\cap \mc O,\,\omega))$ via the injections induced by
  the inclusions of $\mc F\cap \mc O$ into $\mc F$ and into $\mc O$.
  Then the first subset lies in the second iff $\mc F\cap \mc O$ is
  equal to the maximal ideal sheaf $\sM_{\{P\}}$ of some point $P\in C$.
\end{lem}
\begin{proof}
  In any Abelian category, consider two subobjects $A$ and $B$ of an
  object $C$; their intersection is characterized as the kernel of the
  sum map $A\oplus B\to C$.  Thus, as $\mc F$ and $ \mc O$ are
  subsheaves of $\?{\mc O}$, there is a short exact sequence
$$0\to\mc F\cap \mc O\to\mc F\oplus \mc O\to\mc F+ \mc O\to0.$$
To it, apply the left exact functor $H^0(\sHom(\bullet,\,\omega))$.
Thus we obtain the equation
$$H^0(\sHom(\mc F,\,\omega))\bigcap
H^0(\omega)=H^0(\sHom(\mc F+\mc O,\,\omega)),$$ relating the three
subsets of $H^0(\sHom(\mc F\cap \mc O,\,\omega))$.

Therefore, $H^0(\sHom(\mc F,\,\omega))\subset H^0(\omega)$ iff
$H^0(\sHom(\mc F,\,\omega)) = H^0(\sHom(\mc F+\mc O,\,\omega))$.  But
$H^0(\sHom(\mc F,\,\omega))\supset H^0(\sHom(\mc F+\mc O,\,\omega))$.
So to complete the proof, it suffices to prove that $h^0(\sHom(\mc
F,\,\omega))=h^0(\sHom(\mc F+\mc O,\,\omega))$, or $h^1(\mc F)=h^1(\mc
F+\mc O)$.

Consider the following exact cohomology sequence:
$$H^0(\mc F)\to H^0(\mc F+\mc O)\to H^0((\mc F+\mc O)/\mc F)
\to H^1(\mc F)\to H^1(\mc F+\mc O)\to0.$$ By hypothesis, $\mc F$ is a
nonzero proper ideal; so $H^0(\mc F)=0$.  Plainly, $\mc O\subset \mc
F+\mc O\subset\?{\mc O}$, so $h^0(\mc F+\mc O)=1$.  Hence $h^1(\mc
F)=h^1(\mc F+\mc O)$ iff  $h^0((\mc F+\mc O)/\mc F)=1$.  But $(\mc F+\mc
O)/\mc F$ is equal to $\mc O/(\mc F\cap \mc O)$.  And $h^0(\mc O/(\mc
F\cap \mc O))=1$  iff $\mc F\cap \mc O=\sM_{\{P\}}$ for some $P\in C$.
\end{proof}

\begin{thm}[Rosenlicht's Main Theorem]\label{thRMT} Assume $C$ is
  nonhyperelliptic.  Then
  the canonical map is an isomorphism, $\wh\kappa\:\wh C\risom  C'$.
 \end{thm}
\begin{proof}
 Fix $\wh P\in\wh C$, and set $P':=\wh\kappa\wh P$.  Then $\wh\kappa$
provides an inclusion $\smash{\mc O_{C',\,P'}\subset\mc O_{\wh C,\, \wh P}}$,
and we have to prove equality holds.  Now, $\beta\:\wh C\to C$ is an
isomorphism over the Gorenstein locus $G\subset C$ by
Corollary~\ref{coGrlcs}, and $\kappa$ induces an open embedding of $G$
into $C'$ by Theorem~\ref{prGiso}.  Hence, setting $P:=\beta\wh P$, we
may assume $P\notin G$.

For convenience, let $A$ denote the local ring of $P\in C$, and $A'$ the
semilocal ring of $\wh\kappa\beta^{-1}P$.  Lemma~\ref{leR1} provides an
$x\in H^0(\omega)$ with the listed properties (a)--(d).  Property (b)
implies that $\mc O_{\wh C, \wh P}$ is a localization of $A[W_x]$.
Property (c) implies that $A'$ is a localization of $k[W_x]$.
Therefore, to prove the assertion, it suffices to prove that $A\subset
A'$, as then $O_{\wh C, \wh P}\subset O_{C',\,P'}$.  Property (d)
asserts, however, that $A\subset W_x+\sC_P$.  Since $W_x\subset A'$,
therefore
\begin{equation}\label{eqRMT1}
A\subset A'+\sC_P
\end{equation}
Thus it suffices to prove that  $\sC_P\subset A'$.  

Let $\fm\subset A$ be the maximal ideal, and set $\?A:=\?{\mc O}_P$,
which is the integral closure of $A$.  Let's first find an element
$m_0\in A'$ such that $m_0\?A=\fm\?A$.

 Since $\?A$ is a semilocal Dedekind domain, it's a UFD; hence, there
is an $m_1\in\fm$ such that $\?Am_1 = \?A\fm$.  Let
$P_1,\dotsc,P_s\in \?C$ be the points of $\nu^{-1}P$, and
$v_1,\dotsc,v_s$ the corresponding valuations.  Take $n\in \sC_P$ such
that $\?An=\sC_P$.  For $1\le i\le s$, set $a_i:=v_i(n)$ and
$b_i:=v_i(m_1)$; then $a_i\ge b_i\ge1$.  Owing to
Equation~(\ref{eqRMT1}), there exist $m_2\in A'$ and $n_1\in \sC_P$
such that $m_1=m_2+n_1$.  Then, for each $i$, we have $v_i(m_2)\ge b_i$
with equality if $a_i> b_i$.

We have $\delta_P=\eta_P+h^0 (\mc O_P/\sC_P)$ by definition of
$\eta_P$.  But $P\notin G$.  So $\eta_P\ge1$ by Lemma~\ref{leEta}, and $h^0
(\mc O_P/\sC_P)\ge1$ as $P$ is multiple.  Hence $\delta_P\ge2$.

Set $b:=h^0 (\?{\mc O}_P/\sM_{\{P\}})$.  Then $b\ge\delta_P+1$.  So
$b\ge3$.  Moreover,  plainly, $b=\sum b_i$.   

Set $B:=\sum b_iP_i$, which is a divisor on $\?C$.  Fix $i$, and set
$\mc B:=\omega_{\?C}(B-P_i)$.  Then for all $j$, we have
$\deg(P_j+P_i-B)<0$ since $b\ge3$.  So $H^0(\mc O_{\?C}(P_j+P_i-B))=0$.
Hence $H^1(\mc B(-P_j))=0$ by duality.  So the long exact cohomology
sequence becomes
$$0\to H^0(\mc B(-P_j))\to H^0(\mc B)\to k\to 0.$$
So for each $j$, there is an $x_i$ in $H^0(\mc B)$, not in $H^0(\mc
B(-P_j))$.  Since $k$ is infinite, some linear combination of those
$x_i$ is a single $x_i$ that works simultaneously for all $j$.

Set $\mc F:=\nu_*\mc O_{\?C}(P_i-B)$.  Then $\sHom(\mc F,\,\sC\?{\mc
O}\omega)=\nu_*\mc B$.  So $x_i\in H^0(\sHom(\mc F,\,\sC\?{\mc
O}\omega))$.  But $\mc F\supset \nu_*\mc O_{\?C}(-B)$, so $\mc
F\supset\sM_{\{P\}}$.  Plainly, $1\notin\mc F$.  Hence $\mc F\cap\mc
O=\sM_{\{P\}}$.  Therefore, Lemma~\ref{leR2} yields $x_i\in
H^0(\omega)$.

But $x_i\in H^0(\sHom(\mc F,\,\sC\?{\mc O}\omega))$ as $\sC\?{\mc
  O}\omega=\?\omega$ by Lemma~\ref{leCshf}.  And Lemma~\ref{leR1}(a)
implies  $(\?{\mc O}\omega)_P=\?{\mc O}_Px$.  So $x_i\in\Hom(\mc
F_P,\sC_Px)$.  So there is an $f_i\in \Hom(\mc F_P,\,\sC_P)$ such that
$x_i=f_ix$.  Then $f_i\in W_x$ as $x_i\in H^0(\omega)$.  Hence $f_i\in
A'$.

Since $\mc F:=\nu_*\mc O_{\?C}(P_i-B)$, plainly
\begin{equation}\label{eqeqRMT2}
v_j(f_i)\ge\begin{cases}
  a_i-b_i+1& \text{if }j=i,\\
  a_j-b_j& \text{if }j\not=i.
\end{cases}
\end{equation}
In fact, equality holds since $x_i\notin H^0(\mc B(-P_j))$ for all $j$. 

Set $m_3:=f_1^{a_1}\dotsm f_s^{a_s}$.  Then $v_j(m_3)=\sum_i
a_i(a_j-b_j)+a_j$ as equality holds in (\ref{eqeqRMT2}).  Hence
$v_j(m_3)\ge a_j$, with equality iff $a_j=b_j$.

Take $a\in k$, and set $m_0:=m_2+am_3$.  Given $j$, if $a_j>b_j$, then
$v_j(m_2)=b_j$ and $v_j(m_3)\ge a_j$; so then $v_j(m_0)=b_j$ for any
$a$.  If $a_j=b_j$, then $v_j(m_2)\ge b_j$ and $v_j(m_3) = a_j =b_j$; so
then $v_j(m_0)=b_j$ for most $a$.  Thus $m_0\?A=\fm\?A$ as desired.

Second, set $\wh{\sC}:=\sHom(\?{\mc O},\,\wh{\mc O})$ and $\wh A:=\wh{\mc
  O}_P$.  Then $\wh{\sC}_P$ is the conductor of $\?A$ into $\wh A$.
Let's prove that $\sC_P\subsetneq \wh{\sC}_P$.  Indeed,
$\nu_*\omega_{\?C}=\sC\omega$ by Lemma~\ref{leCshf}.  Similarly,
$\wh{\nu}_*\omega_{\?C}=\wh{\sC}\omega_{\wh C}$.  By way of
contradiction, suppose $\sC_P= \wh{\sC}_P$.  Both sides are free
$\?A$-modules as $\?A$ is a UFD.  Hence $\?A\omega_P =
\?A\wh\omega_P$.  Take $f\in\wh\omega_P$ so
that $\?Af=\?A\wh\omega_P$.  Then $f\in \omega_P$ and
$\?Af=\?A\omega_P$.  Hence $\wh Af=\wh A\omega_P$ as $\wh A\omega_P$ is
free.  Therefore,
$$\wh Af\subset \wh\omega_P\subset \omega_P\subset
\wh A\omega_P=\wh Af.$$
 Hence $\wh\omega_P= \omega_P$.  But
$\wh\omega = \sHom(\wh{\mc O},\,\omega)$.  So, by duality,
$\mc O_P=\wh{\mc O}_P$.  Hence $P\in G$, contrary to assumption.  Thus
$\sC_P\subsetneq \wh{\sC}_P$, as desired.

Third, fix $\phi_0 \in \wh{\sC}_P\setminus\sC_P$.  Since $\sC_P$ is an
$\?A$-module, $v_i(\phi_0)<a_i$ for some $i$.  Reordering, we may assume
$i=1$.  Replacing $\phi_0$ by $\phi_0 t$ for a suitable $t\in\?A$, we
may assume
\begin{equation}\label{eqRMT2}
 v_1(\phi_0)=a_1-1 \text{ and }v_j(\phi_0)\ge a_i \text{ for }j\ge 2.
\end{equation}

Fix a $k$-basis $y_1,\dotsc,y_g$ of $H^0(\omega)$, and set
$\phi_j:=y_j/x$ for all $j$.  Then $\phi_1,\dotsc,\phi_g$ form a
$k$-basis of $W_x$.  Now, $\phi_0\in \wh{\sC}_P$, so $\phi_0\in \wh
A$.  Hence $\phi_0$ is a polynomial in $\phi_1,\dots,\phi_g$ with
coefficients in $A$ by Lemma~\ref{leR1}(b).  Say $\phi_0=\sum c_lM_l$
where the $c_l$ belong to $A$ and the $M_l$ are monomials in
$\phi_1,\dots,\phi_g$.    Fix $l$.  Then $M_l\in A'$ by
Lemma~\ref{leR1}(c).  Further, $c_l=c_l'+c_l''$ with $c_l'\in W_x$
and $c_l''\in \sC_P$ by Lemma~\ref{leR1}(d).  So $c_l'\in A'$ by
Lemma~\ref{leR1}(c).  So $c_l'M_l\in A'$.  Further, $c_l''M_l\in\sC_P$
as $\sC_P$ is an $\?A$-module, so an $A'$-module.  Set $\phi_0':=\sum
c_l'M_l$ and $\phi_0'':=\sum c_l''M_l$.  Then $\phi_0'\in A'$ and
$\phi_0''\in \sC_P$.  So $\phi_0'\in \wh{\sC}_P\setminus\sC_P$.  Replace
$\phi_0$ by $\phi_0'$.  Then $\phi_0\in A'$.

Set $\mc G:=\nu_*\mc O_{\?C}(-B)$.  Then $\sHom(\mc G,\,\?\omega)=
\nu_*\omega_{\?C}(B)$.  Plainly, $\mc G\supset\sM_{\{P\}}$ and $1\notin\mc
G$; hence, $\mc G\cap\mc O=\sM_{\{P\}}$.  Therefore, $H^0(\omega_{\?C}(B))
\subset H^0(\omega)$ owing to Lemma~\ref{leR2}. So we may take
$y_1,\dots,y_g$ and $\beta\le g$ so that $y_1,\dots,y_\beta$ belong to
$H^0(\omega_{\?C}(B))$, and yield a basis modulo $H^0(\omega_{\?C})$.

To compute $\beta$, form the long exact cohomology sequence
$$0\to H^0(\omega_{\?C})\to H^0(\omega_{\?C}(B))\to
H^0(\omega_{\?C}(B)/\omega_{\?C})\to H^1(\omega_{\?C})\to
H^1(\omega_{\?C}(B)).$$
Now, $h^1(\omega_{\?C}(B))=0$ because $\deg(B)\ge1$.  Furthermore, plainly,
$h^1(\omega_{\?C})=1$ and $h^0(\omega_{\?C}(B)/\omega_{\?C})=b$.  Hence
$\beta=b-1$.

As noted above, $\?\omega=\sC\?{\mc O}\omega$ by Lemma~\ref{leCshf}.
And $(\?{\mc O}\omega)_P=\?{\mc O}_Px$ by Lemma~\ref{leR1}(a).  Hence
$\phi_i\in \Hom(\mc G_P,\,\sC_P)$ for $1\le i\le \beta$.  So
$v_j(\phi_i)\ge a_j-b_j$ for all $j$.  But $v_j(m_0)=b_j$, so
$v_j(m_0\phi_i)\ge a_j$ for all $j$.  Therefore, $m_0\phi_i\in \sC_P\cap
A'$ for $1\le i\le \beta$.

Recall that $b_1\ge1$.  So (\ref{eqRMT2}) yields $v_j(m_0\phi_0)\ge a_j$
for all $j$.  So $m_0\phi_0\in \sC_P\cap A'$. 

Furthermore, $\phi_0x\notin H^0(\omega_{\?C}(B))$.  Else $\phi_0x$
has no pole outside of the $P_i$.  But equality holds in
(\ref{eqeqRMT2}).  So $\phi_0x$ has a pole of order 1 at $P_1$ and no
other pole.

But no such differential exists.  Indeed, consider the exact sequence
$$0\to H^0(\omega_{\?C})\to H^0(\omega_{\?C}(P_1))\to
H^0(\omega_{\?C}(P_1)/\omega_{\?C})\to H^1(\omega_{\?C})\to
H^1(\omega_{\?C}(P_1)).$$ Here $h^1(\omega_{\?C}(P_1))=0$, and
$h^1(\omega_{\?C})=1$ and $h^0(\omega_{\?C}(P_1)/\omega_{\?C})=1$.
Hence $H^0(\omega_{\?C})$ is equal to $H^0(\omega_{\?C}(P_1))$.  So no
such differential exists.  Thus $\phi_0x\notin H^0(\omega_{\?C}(B))$.

Therefore, $\phi_0x,\phi_1x,\dotsc,\phi_\beta x$ are linearly
independent modulo $H^0(\omega_{\?C})$.  Hence
$\phi_0,\phi_1,\dotsc,\phi_\beta$ are linearly independent modulo
$\sC_P$.  Hence $m_0\phi_0,m_0\phi_1,\dotsc,m_0\phi_\beta$ belong to
$\sC_P\cap A'$, and are linearly independent modulo $\sC_Pm_0$.  But
$\sC_P/\sC_Pm_0$ is isomorphic to $\?A/\?Am_0$, so is of dimension
$b$.  And $b=\beta+1$.  So $m_0\phi_0,m_0\phi_1,\dotsc,m_0\phi_\beta$
yield a basis for $\sC_P\big/\sC_Pm_0$.

For each $n\ge 1$, form the products $m_0^i\phi_j$ for $1\le i\le n$ and
$0\le j\le \beta$.  These products, therefore, yield a basis for
$\sC_P/\sC_Pm_0^n$.  But for $n\gg0$, the ideal $\sC_Pm_0^n$ is
contained in the conductor of $\?A$ into $A'$ since $m_0\?A=\fm\?A$ and
the conductor is an ideal in $\?A$.  But this conductor is also an ideal
in $A'$; so $\sC_Pm_0^n\subset A'$.  But $m_0^i\phi_j\in A'$ for all
$i,j$.  Thus $\sC_P\subset A'$, as desired, and the proof is finally
complete.
\end{proof}

\begin{thm}\label{thNmlty} If $C$ is not Gorenstein, then these six
  conditions are equivalent:
 \begin{enumerate}
 \item $C'$ is projectively normal;
 \item $C'$ is linearly normal;
 \item $C'$ is extremal;
 \item $d'=g'+g-1$;
 \item $C$ is nearly Gorenstein;
 \item $C'=\Spec\bigl(\sHom(\sM_{\{P\}},\,\sM_{\{P\}})\bigr)$ for some point $P$
off the Gorenstein locus.
\end{enumerate}
If \(e)--\(f) hold, then $C$ is of maximal embedding dimension at $P$
iff $C'$ is Gorenstein.  Furthermore, if \(e)--\(f) hold, then $C'$ is
cut out by quadrics and cubics; quadrics suffice if $\eta\ge2$, where
$\eta$ is the invariant of Definition~{\rm\ref{dfCond}}.
 \end{thm}
\begin{proof}
 Conditions (a), (b), and (d) are equivalent by Lemma~\ref{leLnPn}(1),\,(3). 

Suppose (b) holds.  Now, $\wh\kappa$ is an isomorphism by Rosenlicht's Main
Theorem,  Theorem~\ref{thRMT}; also $\wh\kappa^*\mc O_{C'}(1)=\mc O_{\wh C}\omega$.
Hence $h^0(\mc O_{\wh C}\omega)=g$.  So (e) holds by
Lemma~\ref{prNlyG1}.  Conversely, by the same proposition, (e)
implies (b).  Thus (a), (b), (d), and (e) are equivalent.

To prove that (c) and (d) are equivalent, repeat the argument at the end
of the proof of Theorem~\ref{thAn}, mutatis mutandis. 

Given $P\in C$, let $B$ denote the endomorphism ring of its maximal
ideal.  Then $\mc O_P$ is almost Gorenstein, but not Gorenstein, iff
$B=\wh{\mc O}_P$ by Proposition 28 in \cite[p.\,438]{BF}.  And if so,
then $\wh{\mc O}_P$ is also Gorenstein iff $\mc O_P$ is also of maximal
embedding dimension by Proposition 25 in \cite[p.\,436]{BF}.

Since $\beta\: \wh C\to C$ is finite, $\wh C = \Spec(\wh{\mc O})$.
Furthermore, $\beta$ is an isomorphism precisely over the Gorenstein
locus $G\subset C$ by Corollary~\ref{coGrlcs}.  Moreover, $\wh\kappa$ is
an isomorphism by Theorem~\ref{thRMT}.  Hence (e) implies (f) and
the next-to-last assertion.  Conversely, assume (f) holds.  Then $\beta$
is an isomorphism off $P$.  So $P$ is the unique point off $G$.
Furthermore, $\mc O_P$ is almost Gorenstein.  So (e) holds.

To prove the last assertion, note that  $d'=2g'+\eta$ because of (d) and
Proposition~\ref{prdegC'}.  Hence, if $\eta\ge2$, then  $C'$   is cut
out by quadrics owing to Fujita's Corollary~1.14 on p.\,168 in
\cite{Fu83}.

Let $\mc I$ be the ideal of $C'$ in $\bb P^{g-1}$; fix
$l,q\ge0$; and form the long exact sequence
$$H^{q-1}(\mc O_{\bb P^{g-1}}(l))\to H^{q-1}(\mc O_{C'}(l))\to
 H^q(\mc I(l))\to H^q(\mc O_{\bb P^{g-1}}(l)).$$
For $q=1$, the first map is surjective by (a).  For $q=2$, the second
term vanishes for $l\ge1$ by Lemma~\ref{leLnPn}(3) (but not for $l=0$
since $g\ge2)$; for $q\ge3$, this term vanishes as $C'$ is a curve.  For
$q\ge1$, the last term vanishes by Serre's Theorem.  Hence $H^q(\mc
I(3-q))=0$ for $q\ge1$.  Therefore, $C'$ is cut out by quadrics and
cubics by Castelnuovo--Mumford Theory \cite[Prp., p.\,99]{Mu66}.
  The proof is now complete.  \end{proof}

\begin{eg}\label{egCub}
Cubics may be needed to cut out $C'$.  For example, there is a nearly
Gorenstein curve $C$ with $g=4$, with $C'$ smooth, and with $g'=2$; see
Example~\ref{rmCstr}.  Then $d'=5$.  So $C'$ cannot lie on two distinct
quadrics, since their intersection is of degree $4$ by Bezout's Theorem.
\end{eg}

\end{document}